\documentclass[12pt, reqno]{amsart}

\usepackage{amssymb}  
\usepackage{calc}          
\usepackage{multicol}     

\input cyracc.def
\newfam\cyrfam
\font\tencyr=wncyr10 
\font\tenitcyr=wncyi10 

\def\cyr{\fam\cyrfam\twelvecyr\cyracc}
\def\itcyr{\fam\cyrfam\twelveitcyr\cyracc}

\theoremstyle{plain}
\newtheorem{theorem}{Theorem}[section]
\newtheorem{lemma}[theorem]{Lemma}
\newtheorem{definition}[theorem]{Definition}

\newtheorem{corollary}[theorem]{Corollary}
\newtheorem{conjecture}[theorem]{Conjecture}
\theoremstyle{remark}
\newtheorem{example}[theorem]{Example}

\numberwithin{equation}{section}

\newcommand{\seclabel}[1]{\label{sec:#1}} 
\newcommand{\thmlabel}[1]{\label{thm:#1}} 
\newcommand{\lemlabel}[1]{\label{lem:#1}} 
\newcommand{\corolabel}[1]{\label{coro:#1}} 
\newcommand{\deflabel}[1]{\label{def:#1}} 
\newcommand{\exlabel}[1]{\label{ex:#1}} 
\newcommand{\tablelabel}[1]{\label{table:#1}} 
\newcommand{\conjlabel}[1]{\label{conj:#1}} 

\newcommand{\secref}[1]{\ref{sec:#1}} 
\newcommand{\thmref}[1]{\ref{thm:#1}} 
\newcommand{\lemref}[1]{\ref{lem:#1}} 
\newcommand{\cororef}[1]{\ref{coro:#1}} 
\newcommand{\tableref}[1]{\ref{table:#1}} 
\renewcommand{\eqref}[1]{\ref{eq:#1}} 
\newcommand{\conjref}[1]{\ref{conj:#1}} 

\newcommand{\Sym}{\mathrm{Sym}}

\newcommand{\RMlt}{\mathrm{RMlt}}
\newcommand{\RInn}{\mathrm{RInn}}
\newcommand{\LMlt}{\mathrm{LMlt}}
\newcommand{\LInn}{\mathrm{LInn}}
\newcommand{\Mlt}{\mathrm{Mlt}}
\newcommand{\Inn}{\mathrm{Inn}}

\newcommand{\Flex}{\textsc{Flex}}  
  
\newcommand{\RAlt}{\textsc{RAlt}}  
\newcommand{\LAlt}{\textsc{LAlt}}  
\newcommand{\LIP}{\textsc{LIP}}
\newcommand{\RIP}{\textsc{RIP}}
\newcommand{\Mfg}{\textsc{Mfg}}
\newcommand{\PsM}{\textsc{PsM}}
\newcommand{\Ex}{\textsc{Ex}}
\newcommand{\Fex}{\textsc{F}}
\newcommand{\RCC}{\textsc{RCC}}
\newcommand{\LCC}{\textsc{LCC}}
\newcommand{\WIP}{\textsc{WIP}}

\newcommand{\sbl}[1]{\langle#1\rangle}   
\newcommand{\normal}{\trianglelefteq}  
\newcommand{\riv}{^{\rho}} 
\newcommand{\liv}{^{\lambda}} 
\newcommand{\iv}{^{-1}} 
\newcommand{\ld}{\backslash} 
\newcommand{\rd}{/}  

\newcommand{\AAA}{\mathcal{A}}
\newcommand{\LLL}{\mathcal{L}}
\newcommand{\RRR}{\mathcal{R}}

\newcommand{\ZZ}{\mathbb{Z}}

\title[Power-Associative, Conjugacy Closed Loops]
{Power-Associative, Conjugacy Closed Loops}

\author[M.~K.~Kinyon]{Michael~K.~Kinyon}
\address{Department of Mathematical Sciences \\
Indiana University South Bend \\
South Bend, IN 46634 USA}
\email{mkinyon@iusb.edu}
\urladdr{http://mypage.iusb.edu/\symbol{126}mkinyon}
\author[K.~Kunen]{Kenneth~Kunen$^*$}
\thanks{$^*$Partially supported by NSF Grants DMS-0097881 and DMS-0456653}
\address{Department of Mathematics \\
University of Wisconsin \\
Madison, WI 57306 USA}
\email{kunen@math.wisc.edu}
\urladdr{http://www.math.wisc.edu/\symbol{126}kunen}

\date{\today}


\begin{document}

\begin{abstract}
We study conjugacy closed loops (CC-loops) and power-associative
CC-loops (PACC-loops).  If $Q$ is a PACC-loop with nucleus $N$,
then $Q/N$ is an abelian group of exponent $12$; if in addition $Q$
is finite, then $|Q|$ is divisible by $16$ or by $27$.
There are eight nonassociative PACC-loops 
of order $16$, three of which are not extra loops.
There are eight nonassociative PACC-loops of order $27$,
four of which have the automorphic inverse property.

We also study some special elements in loops,
such as Moufang elements, weak inverse property (WIP) elements,
and extra elements.  In a CC-loop, the set of WIP
and the set of extra elements are normal subloops.
For each $c$ in a PACC-loop, $c^3$ is WIP,
$c^6$ is extra, and $c^{12} \in N$.
\end{abstract}

\maketitle

\section{Introduction}
\seclabel{intro}

A loop is \emph{conjugacy closed} (a \emph{CC-loop}) iff it satisfies
the equations:
\[
\begin{array}{clccl}
xy\cdot z = xz\cdot (z\ld (yz)) & (\RCC)
& \qquad &
z\cdot yx = ((zy)\rd z)\cdot zx & (\LCC)
\end{array}
\]
This definition follows Goodaire and Robinson \cite{GRA, GRB};
CC-loops were earlier introduced independently, with different
terminology, by {\cyr So\u\i kis} \cite{SO}. Further discussion
can be found in \cite{DRA1, DRA2, KKP, KUNC}.
The literature is not uniform as to which of these two equations
is left (\LCC) and which is right (\RCC).
With our choice here, also followed in \cite{DRALCC, NS},
\LCC\ is equivalent to saying that the
set of \textit{left} multiplication maps is closed
under conjugation.  In \cite{KKP, KUNC}, 
the equation labels \LCC\ and \RCC\ were arranged
in the opposite order.

In 1982, Goodaire and Robinson \cite{GRA} showed that the nucleus $N(Q)$
of a CC-loop $Q$ is a normal subloop.
A fundamental result in the theory of CC-loops was proved in 1991:

\begin{theorem}[{\cyr Basarab} \cite{BAS}]
\thmlabel{basarab}
Let $Q$ be a CC-loop with nucleus $N = N(Q)$.
Then $Q/N$ is an abelian group.
\end{theorem}

This was conjectured in \cite{GRA}, but
{\cyr Basarab} was apparently unaware of the conjecture,
since he was following the terminology of {\cyr So\u\i kis}.
Because of the differences in terminology, Theorem \thmref{basarab}
was not widely known until recently.
Proofs in English of
Theorem \thmref{basarab} can be found in \cite{DRA1,KKP}.

The CC-loops which are also \emph{diassociative} (that is
every $\sbl{x,y}$ is a group) are precisely the \emph{extra loops}
of Fenyves \cite{FENA,FENB}.
For these, a detailed structure theory was described in \cite{KK}.
The current paper gives a structure theory for those CC-loops
which are only power-associative.

\begin{definition}  For any loop $Q$:
\begin{itemize}
\item[1.]  For $c \in Q$, define $c\riv$ and $c\liv$ by:
$c c\riv = c\liv c = 1$.
\item[2.]  $c \in Q$ is \emph{power-associative}
iff the subloop $\sbl{c}$ is a group.  $Q$ is
\emph{power-associative} iff every element is power-associative.
A \emph{PACC-loop} is a power-associative CC-loop.
\end{itemize}
\end{definition}

The two parts of this definition are related by:

\begin{lemma}[\cite{KUNC}]
\lemlabel{PA-char}
Let $c$ be an element of a CC-loop $Q$. Then
\[
c\text{\ is\ power\ associative}
\qquad \text{iff} \qquad
c\riv = c\liv
\qquad \text{iff} \qquad
c c^2 = c^2 c
\]
\end{lemma}

If $Q$ is a finite nonassociative extra loop, then
$16 \mid |Q|$ (see \cite{KK, KKP}),
and there are exactly five  nonassociative extra loops of order $16$
(see Chein \cite{CH}, p.~49).
We shall show here (Theorem \thmref{16or27}) that if
$Q$ is a finite nonassociative PACC-loop, then 
$16 \mid |Q|$ or $27 \mid |Q|$.
Furthermore (see \S \secref{2loops}),
there are exactly eight nonassociative PACC-loops
of order $16$ (including the five extra loops),
and (see \S \secref{27})
there are exactly eight nonassociative PACC-loops of order $27$.
In \S \secref{ext}, we describe a method of loop extension which will be
useful in constructing the loops in \S\S \secref{2loops},\secref{27}.
In \S \secref{fat}, we prove some results about PACC-loops for which
$Q/N(Q)$ is small; these results will be useful for describing
the PACC-loops of small orders.

In \S\S \secref{WIP},\secref{special}, we discuss some special
kinds of elements in general CC-loops.
Besides the power-associative elements defined above,
there are Moufang and pseudoMoufang elements, extra elements, and WIP
(weak inverse property) elements.
A loop is an extra loop iff all its elements are extra elements;
likewise for Moufang and WIP.
These special elements help us prove some facts about
PACC-loops.  For example (Theorem \thmref{x12}), if $c$ is
in a PACC-loop $Q$, then $c^{12} \in N(Q)$; and this is proved
by showing that $c^3$ is a WIP element and $c^6$ is an extra element.

We shall begin in \S\secref{IMAC} be describing some
basic facts about inner mappings, autotopisms, etc. 

Our investigations were aided by the automated reasoning tool
OTTER \cite{Otter}, and the finite model builder Mace4 \cite{McM},
both developed by McCune. We would also like to thank the referee
for many useful suggestions.

\section{Inner Mappings, Commutators, and Associators}
\seclabel{IMAC}

As usual in a loop $Q$, we define the right and left multiplications
by $xy = x R_y = y L_x$.  These permutations define a number
of important subgroups of $\Sym(Q)$:
\[
\begin{array}{lll}
\Mlt(Q) := \sbl{R_x, L_x : x\in Q} &
\RMlt(Q) := \sbl{R_x : x\in Q} &
\LMlt(Q) := \sbl{L_x : x\in Q} \\
\Inn(Q) := (\Mlt(Q))_1 &
\RInn(Q) := (\RMlt(Q))_1 &
\LInn(Q) := (\LMlt(Q))_1
\end{array}
\]
$\Mlt(Q)$ is called the \emph{multiplication group} of a loop $Q$;
it is generated by all the multiplications.
Then the \emph{inner mapping group} $\Inn(Q)$  is the stabilizer
in $\Mlt(Q)$ of the identity element $1$.
By using only right or only left multiplications, we get
the \emph{right} and \emph{left multiplication groups} and 
the \emph{right} and \emph{left inner mapping groups}.

For $x,y,z\in Q$, define
\[
\begin{array}{c}
	\begin{array}{ccc}
		R(x,y) := R_x R_y R_{xy}\iv & \qquad
		L(x,y) := L_x L_y L_{yx}\iv
	\end{array}
	\\
	T_x := R_x L_x\iv 
\end{array}
\]
These are the standard generators of the inner mapping groups
\cite{Bel, Br}:
\[
\begin{array}{c}
	\Inn(Q) = \sbl{R(x,y), L(x,y), T_x : x,y\in Q} \\
	\begin{array}{ccc}
		\RInn(Q) = \sbl{R(x,y) : x\in Q} & \qquad
		\LInn(Q) = \sbl{L(x,y) : x\in Q}
	\end{array}
\end{array}
\]
Also note that a subloop of $Q$ is normal if and only
if it is invariant under the action of $\Inn(Q)$ \cite{Br}.

\begin{lemma}
\lemlabel{normal}
Let $Q$ be a loop with nucleus $N = N(Q)$ such that
$N \normal Q$ and $Q/N$ is an abelian group. If $M$ is a subloop of $Q$
satisfying $N\leq M$, then $M\normal Q$.
\end{lemma}

\begin{proof}
For each $x\in M$, $\varphi\in \Inn(Q)$, there exists
$n\in N$ such that $(x)\varphi = xn \in M$. Thus
$\Inn(Q)$ leaves $M$ invariant, and so $M$ is normal.
\end{proof}

An \textit{autotopism} of a loop $Q$ is a triple of
permutations $(\alpha,\beta,\gamma)$ such that
$x \alpha \cdot y \beta = (xy)\gamma$ for all $x,y\in Q$.
The autotopisms form a subgroup of $\Sym(Q)^3$.
Defining
\[
\RRR(z) := (R_z,  T_z    ,R_z) \qquad 
\LLL(z) := (  T_z\iv   ,L_z,L_z) \ \ ,
\]
we see that (\RCC) and (\LCC) are equivalent, respectively,
to the assertions that each $\RRR(z)$ and each $\LLL(z)$ is
an autotopism.  In any loop,
if $(\alpha,\beta,\alpha)$ or $(\beta,\alpha,\alpha)$ is an
autotopism and $(1)\alpha = 1$, then $\beta = \alpha$
and $\alpha$ is an automorphism.
Applying this in a CC-loop, where $\RRR(x) \RRR(y) \RRR(xy)\iv$
and $\LLL(x) \LLL(y) \LLL(yx)\iv$ are autotopisms, we get
the following lemma, the parts of which are from
\cite{GRA} and \cite{DRA1}, respectively.

\begin{lemma}
\lemlabel{inner-aut}
For $x,y$ in a CC-loop $Q$,
\begin{enumerate}
	\item[1.] $R(x,y)$ and $L(x,y)$ are automorphisms of $Q$.
	\item[2.] $R(x,y) = T_x T_y T_{xy}\iv$ and $L(x,y) = T_x\iv T_y\iv T_{yx}$,
		so that $\Inn(Q)$ is generated by $\{ T_x : x\in Q\}$.
\end{enumerate}
\end{lemma}

\begin{definition}
Define the \emph{commutator} $[x,y]$
and the \emph{associator} $(x,y,z)$ by:
\[
xy = yx\cdot [x,y] \qquad xy\cdot z = (x\cdot yz)\cdot(x,y,z) \ \ .
\]
The \emph{associator subloop} is $A(Q) := \sbl{(x,y,z) : x,y,z\in Q}$.
\end{definition}

By Theorem \thmref{basarab}, (1) of the following lemma holds for CC-loops.

\begin{lemma}
\lemlabel{T-comm}
In a loop $Q$ with nucleus $N = N(Q)$, the following
are equivalent: 
\begin{enumerate}
	\item[1.] $N\normal Q$ and $Q/N$ is commutative.
	\item[2.] Every commutator is contained in $N$.
\end{enumerate}
In case these conditions hold, $x T_y = x [x,y]$
for all $x,y\in Q$.
\end{lemma}

\begin{proof}
That (1) implies (2) is clear. If (2) holds, then 
$y\cdot (x)T_y = xy = yx [x,y]$, so
$x T_y = x[x,y]$.  Thus, $N T_y \subseteq N$ for all $y$. Since $N$ is
also pointwise fixed by $\RInn(Q)$ and $\LInn(Q)$, we have
$N\normal Q$, and hence (1).
\end{proof}

$A = A(Q)$ is not in general a normal subloop of $Q$; if it is normal,
then $Q/A$ is defined, and it is clearly a group.
The next lemma implies that $A \normal Q$ for CC-loops,
since $A \le N(Q)$ by Theorem \thmref{basarab}.

\begin{lemma}
\lemlabel{assoc-normal}
Let $Q$ be a loop with associator subloop $A = A(Q)$.
If $A\leq N(Q)$, then for all $x,y,z,u\in Q$,
\begin{enumerate}
	\item[1.] $(x,y,z) T_u = (x,yz,u)\iv (y,z,u)\iv (x,y,zu) (xy,z,u)$
	\item[2.] $[(x,y,z),u] = (x,y,z)\iv (x,yz,u)\iv (y,z,u)\iv (x,y,zu) (xy,z,u)$
\end{enumerate}
In particular, $A(Q) \normal Q$.
\end{lemma}

\begin{proof}
Compute
\begin{align*}
	(x\cdot yz) (x,y,z)u
		&= (xy\cdot z)u
		= xy\cdot zu\cdot (xy,z,u)
		= x(y\cdot zu)(x,y,zu)(xy,z,u) \\
		&= x(yz\cdot u)(y,z,u)\iv (x,y,zu)(xy,z,u) \\
		&= (x\cdot yz)\cdot u\cdot (x,yz,u)\iv (y,z,u)\iv (x,y,zu)(xy,z,u).
\end{align*}
Cancel $x\cdot yz$ and then divide both sides on the left by $u$
to obtain (1). (2) follows from (1), using $(x,y,z) \in N$.
By (1), $a T_u \in A$ for each generator $a$ of $A$ and each $u \in Q$.
It follows that each $A T_u \subseteq A$; to verify this,
use $A\leq N(Q)$ and note that 
$(mn) T_u = m T_u \cdot n T_u$ whenever $m,n, m T_u, n T_u \in N$.
Thus, $A$ is normal, since any subloop
of $N(Q)$ is pointwise fixed by $\RInn(Q)$ and $\LInn(Q)$. 
\end{proof}

\begin{lemma}
\lemlabel{Q/N-group}
Let $Q$ be a loop with nucleus $N = N(Q)$. If $N\normal Q$ and 
$Q/N$ is a group, then
\begin{itemize}
	\item[1.] $(x,y,z) = (ux,vy,wz)$ for all $x,y,z \in Q$ and $u,v,w \in N$.
	\item[2.] $A(Q) \le Z(N)$ so that $A(Q)$ is an abelian group.
	\item[3.] $(x^\rho, y, z) = (x^\lambda, y, z)$.
\end{itemize}
\end{lemma}

\begin{proof}
For (1) and (2), see, for instance, \cite{KKP}, \S5.
For (3), note that $x\riv \rd x\liv \in N$.
\end{proof}

In particular, Lemma \lemref{Q/N-group} applies to CC-loops
by Theorem \thmref{basarab}.

\begin{lemma}
\lemlabel{assoc-CC}
In a loop,
\begin{enumerate}
\item[1.] \RCC\ is equivalent to $(x,y,z) = (x,z, yT_z)$.
\item[2.] \LCC\ is equivalent to $(x,y,z) = (yT_x\iv,x,z)$.
\end{enumerate}
\end{lemma}

\begin{proof}
$\RCC$ is equivalent to $xy \cdot z =  xz\cdot yT_z$, so (1) is clear from:
\begin{align*}
 xy\cdot z =& (x\cdot yz)\cdot(x,y,z) \\
 xz\cdot yT_z =& (x\cdot (z \cdot yT_z))\cdot(x,z,yT_z) 
= (x\cdot yz)\cdot (x,z,yT_z)
\end{align*}
(2) is the mirror of (1).
\end{proof}

\begin{lemma}
\lemlabel{ccequiv}
A loop $Q$ is a CC-loop iff
\begin{itemize}
\item[a.] All associators are invariant under permutations of their arguments,
and
\item[b.] All commutators are nuclear.
\end{itemize}
\end{lemma}
\begin{proof}
For $\Rightarrow$, (a), which is (\cite{KKP},
Theorem 4.4), follows from Theorem \thmref{basarab} and
Lemmas  \lemref{Q/N-group} and
\lemref{assoc-CC} (since $yT_z = vy$ for some $v \in N$),
while (b) follows from Theorem \thmref{basarab}.
For $\Leftarrow$, note that (a), (b), and Lemmas \lemref{T-comm}
and \lemref{Q/N-group} yield
$(x, y, z) = (x, z, yT_z)$ and
$(x, y ,z) = (yT_x\iv, x, z)$, and then apply
Lemma \lemref{assoc-CC}.
\end{proof}

We shall frequently use Lemmas
\lemref{Q/N-group} and \lemref{ccequiv}
without comment when writing equations with associators.
We next collect some further properties of associators
in CC-loops. 

\begin{lemma}
\lemlabel{assoc-props}
In any CC-loop $Q$:
\begin{itemize}
	\item[1.] $(xy,z,u) = (x,z,u)T_y\cdot (y,z,u) = (x,z,u)\cdot (y,z,u)T_x$
	\item[2.] $[(x,z,u),y] = [(y,z,u),x]$
	\item[3.] $(x\riv,z,u) T_x = (x,z,u)\iv$
	\item[4.] $(xy,z,u) = 1$ iff $(x,z,u) = (y\riv,z,u)$
	\item[5.] for each $x,y\in Q$, $\{u : (u,x,y) = 1\}$ is a subloop of $Q$.
\end{itemize}
\end{lemma}

\begin{proof}
We know that associators are invariant under permutations of their arguments
(Lemma \lemref{ccequiv}) and lie in the center of the nucleus
(Lemma \lemref{Q/N-group}).  Also, 
$(xy,z,u) = (yx,z,u)$ because commutators are nuclear.  Now,
by Lemma \lemref{assoc-normal}(1),
\begin{align*}
(x,y,z) T_u = (x,yz,u)\iv (y,z,u)\iv (x,y,zu) (xy,z,u) \\
(x,z,y) T_u = (x,zy,u)\iv (z,y,u)\iv (x,z,yu) (xz,y,u)\,.
\end{align*}
But these are equal, so
\[
(x,y,zu) (xy,z,u) = (x,z,yu) (xz,y,u)\,.
\]
Using this and Lemma \lemref{assoc-normal}(1) again, we get (1):
\begin{align*}
(x,z,u) T_y \cdot  (y,z,u)  = &(x,zu,y)\iv  (x,z,uy) (xz,u,y) = \\
 &(x,zu,y)\iv  (x,y,uz) (xy,u,z) = (xy,z,u)\,.
\end{align*}
(2) follows from Lemma \lemref{T-comm} and (1).
We obtain (3) by taking $y = x\riv$ in (1). 
(4) follows from (1) and (3).
In (5), note that $\{u : (u,x,y) = 1\} = \{u : u R(x,y) = u\}$,
which is a subloop because $R(x,y)$ is an automorphism.
\end{proof}

\begin{lemma}
\lemlabel{inner-assoc}
For all $x,y,z$ in a CC-loop $Q$,
\begin{enumerate}
	\item[1.] $z L(x,y) = z (z,x,y)\iv$
	\item[2.] $z R(x,y) = z (z,x\liv,y\liv)$
	\item[3.] $R(x,y)\iv = L(x\liv,y\liv)$,
		so that $\RInn(Q) = \LInn(Q)$.
	\item[4.] $R(x,y) = R(y,x)$
	\item[5.] $R(x,y) R(u,v) = R(u,v) R(x,y)$,
		so that $\RInn(Q)$ is an abelian group.
\end{enumerate}
\end{lemma}

\begin{proof}
(1) is from (\cite{KKP}, \S4). Next, we compute
\[
x R(y,z) R_{yz} = xy\cdot z = x R_{yz}\cdot(x,y,z)
\]
so that $x R(y,z) = x R_{yz} R_{(x,y,z)} R_{yz}\iv
= x [(yz \cdot (x,y,z)) \rd (yz)]$. By
Lemma \lemref{assoc-props}(4),
$yz \cdot (x,y,z) = y (x,y,z\liv)\iv z
= (x,y\liv,z\liv)\cdot yz$, and (2) holds.
(3) follows from (1), (2), and Lemma \lemref{inner-aut}(1).
(4), which is from \cite{KKP}, follows from (3) and
Lemma \lemref{assoc-props}(1). (5), which is also
from \cite{KKP}, follows from (2) and Lemmas
\lemref{inner-aut}(1) and \lemref{Q/N-group}(1). 
\end{proof}

By finding an expression for $L(x,y)$ as a product of right
multiplications, Dr\'{a}pal \cite{DRA1} was the first to show
that $\RInn(Q) = \LInn(Q)$ for a CC-loop $Q$. However, the
equation in Lemma \lemref{inner-assoc}(3) relating the
generators of the two groups seems to be new.

The following inner mapping was used also in \cite{KKP, KUNC}:
\[
E_x := R(x,x\riv) = R_x R_{x\riv} \ \ .
\]
The next lemma collects some of its properties:

\begin{lemma}
\lemlabel{Ebasic}
For every $x,y$ in a CC-loop $Q$:
\begin{itemize}
	\item[1.] Each $E_x$ is an automorphism.
	\item[2.] $y E_x = y (y,x\liv,x)$.
	\item[3.] $E_x = R(x\liv,x) = L(x,x\liv)\iv = L(x\riv,x)\iv = R_x L_x R_x\iv L_x\iv$ .
	\item[4.] If $x$ is power-associative, then $E_{x^n} = E_x^{n^2}$,
$[E_x, L_x ] = [E_x, R_x] = I$, and
$R_{x^n} = R^n_x E_x^{(n-1)n/2}$,
$L_{x^n} = L^n_x E_x^{-(n-1)n/2}$,
$R_x^n L_x^m = L_x^m R_x^n E_x^{mn}$ .
	\item[5.] If $Q$ is a PACC-loop, then $E_x^6 = I$.
\end{itemize}
\end{lemma}

\begin{proof}
(1) and (2) are just specializations of Lemmas
\lemref{inner-aut}(1) and \lemref{inner-assoc}(2),
respectively.
The first two equalities of (3) follow from (2) and Lemma
\lemref{inner-assoc}(1)(2); that is,
$y R(x\liv,x) = y (y,x,x\liv)$ and
$y L(x,x\liv)\iv = y (y,x,x\liv)$. The third equality
follows from these and Lemma \lemref{Q/N-group}(3).
The remaining equality
of (3) is from \cite{KUNC}, as is (4). (5) is from \cite{KKP}.
\end{proof}

We conclude this section with the following easy criterion
for checking that a subset is a subloop:

\begin{lemma}
\lemlabel{subloop}
A subset $X$ of a CC-loop $Q$ is a subloop iff $X$ is closed under
product and either $\rho$ or $\lambda$.
\end{lemma}

\begin{proof}
By Lemma \lemref{Ebasic}, $L_x\iv = E_x L_{x\riv} = L_{x\liv} E_x$
and $R_x\iv = R_{x\riv} E_x\iv = E_x\iv R_{x\liv}$. 
This yields the equations
\[
\begin{array}{lcc}
x \ld y &= x\riv (yx\cdot x\riv) &= (x\liv y\cdot x\liv) x \\
y \rd x &= x(x\riv\cdot y x\riv) &= (x\liv \cdot xy) x\liv
\end{array}
\]
which immediately yield the desired result. 
\end{proof}

\section{WIP Elements}
\seclabel{WIP}

The role played by weak inverse property elements in
CC-loops was already highlighted in \cite{GRB, KKP}. In
this section we elaborate further on that theme.

\begin{definition}
\deflabel{WIP}
An element $c$ of a loop $Q$ is a
\emph{weak inverse property (WIP) element} iff 
for all $x\in Q$,
\[
c(xc)\riv = x\riv \qquad\qquad 
(c x)\liv c = x\liv\,.
\tag{\WIP}
\]
Let $W(Q)$ denote the set of all WIP elements of Q.
\end{definition}

\noindent The two equations defining a WIP element are equivalent
in all loops (\cite{KKP}, Lemma 2.18). Also note that
$N(Q)\subseteq W(Q)$.

\begin{lemma}
\lemlabel{wip}
For an element $c$ of a CC-loop, the following are equivalent;
in $(ii)$--$(vii)$, the variable $x$ is understood to be universally quantified:
\begin{center}
\begin{tabular}{c}
i. \quad $c$ is a WIP element \\
\begin{tabular}{lcclc}
ii. & $x (cx)\riv = c\riv$ & \qquad &
iii. & $(xc)\liv x = c\liv$ \\
iv. & $c = (c \cdot xE_c) x\riv$ & \qquad &
v. & $x = (x \cdot cE_x) c\riv$ \\
vi. & $(c,x,x\riv) = (x\riv,c,c\riv)$ & \qquad &
vii. & $ (x,c,c\riv) = (c\riv,x,x\riv)$.
\end{tabular}
\end{tabular}
\end{center}
\end{lemma}

\begin{proof}
(i) holds iff $x\cdot c(xc)\riv = 1$, that is, iff
$c\cdot [c\ld (xc)] (xc)\riv = 1$ for all $x$, using \LCC.
Replacing $x$ with $(cx)\rd c$, we have that (i) holds iff
$c\cdot x(cx)\riv = 1$ for all $x$. Thus $(i)\leftrightarrow (ii)$,
and the mirror of this argument yields $(i) \leftrightarrow (iii)$.

Next, (i) holds iff $(xc)\riv = c\ld x\riv$; that is, $1 = xc\cdot (c\ld x\riv)$.
Multiplying on the left by $c$ and using \LCC, we have that (i) holds iff
$c = xR_c L_c R_c\iv \cdot x\riv = (c\cdot xE_c)\cdot x\riv$, recalling
$E_c = R_c L_c R_c\iv L_c\iv$ (see Lemma \lemref{Ebasic}).
Thus $(i)\leftrightarrow (iv)$. Interchanging $c$
and $x$ in this argument, we get $(ii) \leftrightarrow (v)$.

For $(iv) \leftrightarrow (vi)$, use Lemma \lemref{Ebasic}(2)
plus Lemma \lemref{assoc-props}(3) to get:
\[
(c\cdot xE_c)x\riv = 
cx \cdot (x, c\riv, c) \cdot x\riv =
cx \cdot x\riv \cdot (x\riv, c\riv, c)\iv =
c \cdot (c,x,x\riv) (x\riv, c\riv, c)\iv \ \ .
\]
Interchanging $c$ and $x$ in this argument yields $(v) \leftrightarrow (vii)$.
\end{proof}

\begin{corollary}
\corolabel{wip2}
For a WIP element $c$ of a CC-loop, $x E_c = x$ iff $c E_x = c$.
\end{corollary}

\begin{proof}
By parts $(iv)$ and $(v)$ of Lemma \lemref{wip}.
\end{proof}

\begin{theorem}
\thmlabel{W-subloop}
In a CC-loop $Q$, $W(Q)$ is a normal subloop.
\end{theorem}

\begin{proof}
We show that $W = W(Q)$ is a subloop, so fix $b, c\in W$,
and we show that $W$ contains $bc$ and $c\riv$ (see Lemma \lemref{subloop}).
Normality will follow from Lemma \lemref{normal}.

For $bc$: Set $u = c\cdot xc\riv$, and note that
$bu\cdot c = bc\cdot x$ by \RCC, since $R_{c\riv} L_c R_c L_c\iv = I$
(by Lemma \lemref{Ebasic}(3):
$E_c = R_c R_{c\riv} =  R_c L_c R_c\iv L_c\iv$ ).
Then, using 
Lemma \lemref{wip},
$ R_{c\riv} L_c = L_c R_c\iv$,
\LCC, and $c \in W$, we have
\[
b\riv = u(bu)\riv =
((cx)\rd c)\cdot  (bu)\riv  =
c \cdot x (c\ld (bu)\riv )
= c\cdot x(bu\cdot c)\riv = c\cdot x(bc \cdot x)\riv .
\]
Now
\[
bc = (b\riv)\liv c = (c\cdot x(bc \cdot x)\riv)\liv c = (x(bc \cdot x)\riv)\liv ,
\]
that is, $(bc)\riv = x(bc\cdot x)\riv$. By Lemma \lemref{wip}, $bc\in W$.

For $c\riv$: By \LCC, 
\[
c [(c\ld x) x\riv] = (x\rd c)\cdot cx\riv 
= (cx\riv)\liv \cdot cx\riv = 1 \ \ .
\]
Then, using Lemma \lemref{Ebasic}(3)
($E_{c\riv} = L({c\riv},c)\iv = L_c\iv L_{c\riv}\iv  $):
\[
c\riv = (c\ld x) x\riv =
x E_{c\riv} L_{c\riv} \cdot x\riv = (c\riv\cdot xE_{c\riv}) x\riv \ \ ,
\]
so $c\riv$ is WIP by Lemma \lemref{wip}.
\end{proof}

\begin{corollary}
\corolabel{wip-cyclic}
In a CC-loop $Q$, if $c\in W(Q)$, then $\sbl{c}\subseteq W(Q)$.
\end{corollary}

For a PACC-loop, we can say something about the structure of the
quotient by the WIP subloop. To this end, we quote the
following from \cite{KKP}, Theorem 8.4:

\begin{theorem}
\thmlabel{cubes}
Let $Q$ be a PACC-loop. For each $c\in Q$, $c^3 \in W(Q)$.
\end{theorem}

\begin{corollary}
\corolabel{Q/W}
Let $Q$ be a PACC-loop with WIP subloop $W(Q)$. Then $Q/W(Q)$
is an elementary abelian $3$-group.
\end{corollary}

Then, since $N\leq W$, we have:

\begin{corollary}
\corolabel{nucwip}
If $Q$ is a PACC-loop and 
$a^r\in N(Q)$, where $\gcd(r,3) = 1$, then $a$ is a WIP element.
\end{corollary}

WIP elements have the following further associator properties:

\begin{lemma}
\lemlabel{wip-assoc}
In a CC-loop, if $a$ is a WIP element and $b$ is arbitrary,
then $(a^2,x,y) = (b^2,u,v) = 1$ for all $x,y \in \sbl{a,b}$
and all $u,v\in\sbl{a}$.  If $b$ is also power-associative,
then $(b^2,u,v) = 1$ for all $u,v \in \sbl{a,b}$.
\end{lemma}

\begin{proof}
By Lemma \lemref{wip}\emph{(vi)} and \emph{(vii)}, we have
\[
 (a,b,b\riv) = (b\riv,a,a\riv)\qquad (a\riv,b\riv,b) = (b,a,a\riv) 
\]
and
\[
  (b,a,a\riv) = (a\riv,b,b\riv) \qquad (b\riv,a,a\riv) = (a\riv,b,b\riv)
\]
so all these associators are equal.  It follows by 
Lemma \lemref{assoc-props}(4) that
$(a^2,b,b\riv) = 1$ and $(b^2, a, a\riv) = 1$.
By Lemma \lemref{assoc-props}(5), we have
$(a^2,x,y) = (b^2,u,v) = 1$ for all $x,y \in \sbl{b}$
and all $u,v\in\sbl{a}$.  

Since $b$ was arbitrary, we can also replace $b$ by $a$ to get
$(a^2,x,y) = 1$ for all $x,y \in \sbl{a}$.  Then, by
Lemma \lemref{assoc-props}(1),
\[
(a^2,a\riv,b) = (a,a\riv,b)T_a\cdot (a,a\riv,b) =
(a,b\riv,b)T_a\cdot (a,b\riv,b) =
(a^2,b\riv,b) = 1 \ \ ,
\]
so by Lemma \lemref{assoc-props}(5), we have
$(a^2,x,y) = 1$ for all $x\in \sbl{a}$ and $y \in \sbl{b}$.
Applying Lemma \lemref{assoc-props}(5) two more times, we get
$(a^2,x,y) = 1$  for all $x,y \in \sbl{b}$.

If $b$ is also power-associative,
then $(b^2,u,v) = 1$ for all $u,v \in \sbl{b}$, and the above
argument then gives us  $(b^2,u,v) = 1$ for all $u,v \in \sbl{a,b}$.
\end{proof}

This implies the following properties of $2$-generated
CC-loops when one of the generators has WIP. Parts of
this lemma are in \cite{KKP} (see Theorems 7.8, 7.10);
the proof here is different.

\begin{lemma}
\lemlabel{wip-2gen}
Let $Q$ be a CC-loop with nucleus $N = N(Q)$, and assume that
$Q = \sbl{a,b}N$ where $a$ has WIP. Then $a^2\in N$. 
Further,
\begin{enumerate}
\item[1.] If $b$ is power-associative, then $b^2\in N$
and $\sbl{a^2,b}$ is a group.
\item[2.] If $a$ and $b$ are power-associative, then
$\sbl{a,b^2}$ is a group,
$(a,a,b)=(a,b,b)$ generates $A(Q)$,
$|A(Q)| \le 2$, $Q$ is a PACC-loop, and $A(Q)\leq Z(Q)$.
\end{enumerate}
\end{lemma}

\begin{proof}
First, $a^2 \in N$ by Lemma \lemref{wip-assoc} because
$(a^2,x,y) = 1$ for all $x,y \in Q$.
Likewise, if $b$ is power-associative then 
$b^2 \in N$ by Lemma \lemref{wip-assoc},
and $\sbl{a^2, b} \le \sbl{\{b\} \cup N} = \sbl{b}N$ is a group.

Now assume that both $a$ and $b$ are power-associative.
Then $\sbl{a,b^2}$ is a group because $b^2\in N$.
By Theorem
\thmref{basarab} and Lemmas \lemref{T-comm} and
\lemref{Q/N-group}, every associator in $Q$ can be
expressed in the form $(a^i b^j, a^k b^{\ell}, a^m b^n)$
for integers $i,j,k,{\ell},m,n$.
Furthermore, since $a^2, b^2\in N(Q)$,
Lemma \lemref{Q/N-group} implies
the value of  $(a^i b^j, a^k b^{\ell}, a^m b^n)$
depends only on $i,j,k,{\ell},m,n \bmod 2$.

Lemma \lemref{wip} $(vi)$ and $(vii)$ now gives us
$(a,a,b) = (a,b,b)$; call this value $u$.  Then
\[
u = (a,a,a)T_b \cdot (b,a,a) =
(ab,a,a) 
= (a,a,a) \cdot (b,a,a) T_a = u T_a
\]
by Lemma \lemref{assoc-props}(1), so $(ab,a,a) = u T_a = u$;
similarly, $(ab,b,b) = u T_b = u$. Then,
\[
u^2 =  (a,a,b) \cdot (a,a,b) T_a = (a^2,a,b) = 1 \ \ .
\]
Then
$ (ab,a,b) = (a,a,b) \cdot (b,a,b) T_a = u^2  = 1$ and
$ (ab,a,ab) = (a,a,ab) \cdot (b,a,ab) T_a = u$, and likewise
$ (ab,a,ab) = u$.
Finally, 
$ (ab,ab,ab) = (a,ab,ab) \cdot (b,ab,ab) T_a = u^2 = 1$.
This accounts for all associators, so $A(Q) = \{1,u\}$.

Since $(x,x,x) = 1$ for each of $x = a,b,ab$,
we have $(x,x,x) = 1$ for all $x$, so $Q$ is a PACC-loop.
Finally $u \in Z(Q)$ because $u$ commutes with $a$ and $b$.
\end{proof}

\begin{corollary}[\cite{KKP}, Corollary 8.5]
\corolabel{old-groups}
For each $a,b$ in a PACC-loop, $\sbl{a^3,b^2}$
and $\sbl{a^6,b}$ are groups.
\end{corollary}

\begin{proof}
By Theorem \thmref{cubes}, $a^3$ has WIP, so 
the result follows from Lemma \lemref{wip-2gen}.
\end{proof}

\section{Moufang, PseudoMoufang, and Extra Elements}
\seclabel{special}

\begin{lemma}
\lemlabel{eqns}
Let $Q$ be a CC-loop, and fix $a,b\in Q$.
If any of the following hold, then they all hold: 
\[
\begin{array}{c}
a\cdot ba = ab\cdot a \qquad {\rm (\Flex)} \\
\begin{array}{ccccc}
a\cdot ab = a^2 b & {\rm (\LAlt)} && ba\cdot a = ba^2 & {\rm (\RAlt)} \\
a\liv \cdot ab = b & {\rm (\LIP)} &\qquad& ba\cdot a\riv = b & {\rm (\RIP)} \\  
\forall x \,[ ab\cdot xa = a(bx\cdot a)]]  & {\rm ({\Mfg}1)} && 
        \forall x \,[ ax\cdot ba = (a\cdot xb)a ] & {\rm ({\Mfg}2)} \\
\forall x \,[ ab\cdot ax = a(ba\cdot x)]  & {\rm ({\Fex}1)} && 
        \forall x \,[ xa\cdot ba = (x\cdot ab)a ] & {\rm ({\Fex}2)}
\end{array}
\end{array}
\]
\end{lemma}

\begin{proof}
(\Flex), (\LAlt), and (\RAlt) are equivalent by 
Lemma \lemref{ccequiv}. (\RIP) is $bE_a = b$, and
Lemma \lemref{Ebasic}(3) shows
that this is equivalent to (\LIP) and to (\Flex).
Taking $x = 1$ in ({\Mfg}1) or in ({\Fex}1) gives (\Flex).
To see that $(\Flex) \to (\Fex1 \wedge {\Mfg}1)$,
apply \LCC\ and (\Flex) to get:
\[
a (ba \cdot a\ld u) = ((a\cdot ba )\rd a)\cdot (a \cdot a\ld u) =
ab \cdot u \ \ .
\]
Set $u = ax$ to get ({\Fex}1).  Setting $u = xa$ and applying \RCC\ yields
({\Mfg}1) by
\[
ab \cdot xa = a (ba \cdot a\ld xa) =   a (bx\cdot a ) \ \ .
\]
The equivalence with ({\Mfg}2) and with ({\Fex}2) is established similarly.
\end{proof}

If $V$ is a variety of loops defined by some (universally quantified)
equation, we may fix one variable of the equation and define an element
$c$ of an arbitrary loop to be a ``$V$ element'' iff $c$ satisfies
that equation with the other variables quantified. For short, we
simply say ``$c$ is $V$". We have already seen one example of this
with WIP elements. For another, one of the Moufang equations is 
$\forall x,y,z \,[(z \cdot  xy) z = zx\cdot yz]$, and
Bruck  defines $c$  to be a Moufang element iff
$\forall x,y \,[(c \cdot  xy) c = cx\cdot yc]$
(see \cite{Br},VII\S2).
Of course, $V$ may be defined by several equivalent equations,
each of which may have several variables, so the proper definition of
a ``$V$ element'' is not uniquely determined.
The following definitions seem to be useful in CC-loops.

\begin{definition}
\deflabel{elements}
An element $c$ of a CC-loop $Q$ is a\textup{(}n\textup{)}
\begin{center}
\begin{tabular}{rclc}
\emph{Moufang element} & \quad iff \quad & 
$\begin{cases}
  \forall x,y \, [ c(xy\cdot c) = cx\cdot yc] \\
   \forall x,y \, [(c\cdot xy)c = cx\cdot yc]
  \end{cases}$ 
& {\rm (\Mfg)}
\\
\emph{pseudoMoufang element} &iff&
$\begin{cases}
  \forall z,x \, [z(cx\cdot z) = zc\cdot xz] \\
  \forall z,x \, [(z\cdot xc)z = zx\cdot cz] 
  \end{cases}$
& {\rm (\PsM)}
\\
\emph{extra element} &iff&
$\forall x,y \, [c(x\cdot yc) = (cx\cdot y)c]$
& {\rm (\Ex)}
\end{tabular}
\end{center}
Let $M(Q)$, $P(Q)$, and $Ex(Q)$ denote the
sets of Moufang, pseudoMoufang, and extra elements of $Q$,
respectively.
\end{definition}

\noindent
In all loops, the two equations defining a Moufang element
are equivalent because each one
implies $\forall z \, [c \cdot zc = cz \cdot c]$.
Note that $(\PsM)$ is obtained by fixing a different variable in
the Moufang laws.  In CC-loops,
the two equations defining a Pseudo\-Moufang element
are equivalent because, by Lemma \lemref{eqns},
$(\Mfg 1) \leftrightarrow (\Mfg 2)$.

Note that elements of $M(Q)$, $P(Q)$, and $Ex(Q)$
are power-associative.

The following is immediate from Lemmas \lemref{eqns}
and \lemref{Ebasic}(2).

\begin{corollary}
\corolabel{E-Mfg}
Let $c$ be an element of a CC-loop $Q$.
\begin{center}
\begin{tabular}{lccccc}
$(i)$ & $c \in M(Q)$ & iff & $E_c = I$ & iff &
$(x,c,c\riv) = 1$ for all $x$. \\
$(ii)$ & $c \in P(Q)$  & iff & $cE_x = c$ for all $x$ &
iff & $(c,x,x\riv) = 1$ for all $x$.
\end{tabular}
\end{center}
\end{corollary}

This corollary plus Lemma \lemref{assoc-props}(5)
yields the next two corollaries:

\begin{corollary}
\corolabel{mfg-group}
In a CC-loop $Q$, if $a \in M(Q) \cap P(Q)$
and $b$ is a power-associative element, then $\sbl{a,b}$
is a group.
\end{corollary}

\begin{corollary}
\corolabel{mfg-more}
In a CC-loop $Q$, if $a \in M(Q)$ then $\sbl{a} N(Q) \subseteq M(Q)$.
\end{corollary}

\begin{theorem}
\thmlabel{other-sblp}
In a CC-loop $Q$, $P(Q)$ is a normal subloop.
\end{theorem}

\begin{proof}
$P(Q)$ is exactly the fixed point subset of the
automorphisms $E_x$ (Corollary \cororef{E-Mfg}\emph{(ii)}), and thus is a subloop. Since
$P(Q) \geq N(Q)$, the normality follows from
Lemma \lemref{normal}.
\end{proof}

\begin{lemma}
\lemlabel{M-products}
Let $a,b$ be Moufang elements of a 
CC-loop $Q$. Then the following are equivalent:
\begin{center}
\begin{tabular}{lcclc}
1. & $ab$ is a Moufang element &&
2. & $ba$ is a Moufang element \\
3. & $(x,a,b) = (x,a\riv,b)$ for all $x\in Q$ &&
4. & $(x,a,b) = (x,a,b\riv)$ for all $x\in Q$
\end{tabular}
\end{center}
\end{lemma}

\begin{proof}
$(1) \leftrightarrow (2)$ holds by Corollary \cororef{mfg-more},
because $[a,b] \in N(Q)$.  For $(1) \leftrightarrow (3)$,
note that $ab \in M(Q)$ iff $(x,ab,a\riv b\riv) = 1$ for all $x$
(by Corollary \cororef{E-Mfg}), and
\begin{align*}
&(x,ab,a\riv b\riv) =
(x,a,a\riv b\riv)T_b \cdot (x,b,a\riv b\riv) = \\
&(x,a,a\riv ) T_{b\riv} T_b \cdot (x,a, b\riv)T_b
\cdot (x,b,b\riv)T_{a\riv} \cdot (x,b,a\riv) =
  (x,a, b)\iv \cdot (x,b,a\riv)  \ \ ,
\end{align*}
using $a,b \in M(Q)$ and Lemma \lemref{assoc-props}.
Similarly, $(2) \leftrightarrow (4)$.
\end{proof}

While we will use Lemma \lemref{M-products} below, we have not been
able to show that $M(Q)$ is a subloop even in the power-associative
case. Thus we offer the following.

\begin{conjecture}
\conjlabel{moufang}
There exists a PACC-loop $Q$ in which $M(Q)$ is not a subloop.
\end{conjecture}

Note that by Lemma \lemref{wip-extra} below, such a $Q$ cannot
be WIP.

Extra loops, introduced by Fenyves \cite{FENA,FENB}, are loops satisfying
one of the following three equivalent identities:
\[
 z  x \cdot z  y = z  (x  z \cdot y)  \qquad
 x  z \cdot y  z= (x \cdot z  y)  z  \qquad
 z (x\cdot yz) = (zx\cdot y)  z \ \ .
\]
The identities ({\Fex}1) and ({\Fex}2) in Lemma \lemref{eqns} are obtained
by fixing variables in the first two of these,
but for our definition of ``extra element'', we found it more useful
to fix a variable in the third one.
Extra loops are both CC and Moufang, and all squares
in an extra loop are nuclear.  Generalizing this,

\begin{lemma}
\lemlabel{nuc-square}
For an element $c$ of a CC-loop $Q$, the following are equivalent:
\begin{center}
\begin{tabular}{lclclc}
& & $(i)$ & $c$ is extra & & \\
$(ii)$ & $\forall x, y \,[c(x\cdot cy) = (cx\cdot c)y]$ & & &
       $(iii)$ & $\forall x, y \,[(yc\cdot x)c = y(c\cdot xc)]$ \\
$(iv)$ & $\forall x, y \,[c(x\cdot cy) = (c\cdot xc)y]$ & & &
       $(v)$ & $\forall x, y \,[(yc\cdot x)c = y(cx\cdot c)]$
\end{tabular}

$(vi)\quad c$ is Moufang and $c^2 \in N(Q)$
\end{center}
\end{lemma}
\begin{proof}
First note that (i)--(v) imply that $c \cdot xc = cx \cdot c$
(setting $y = 1$ in (\Ex) or (ii)--(v)).
Thus, by Lemma \lemref{eqns},
in all cases, $c$ is Moufang, and hence also 
$c\cdot cx = c^2 x$ (i.e., $L_c^2 = L_{c^2}$).
(ii)$\leftrightarrow$(iv) and
(iii)$\leftrightarrow$(v) follow from $c \cdot xc = cx \cdot c$.

For (i)$\leftrightarrow$(vi), note that 
$c$ is extra iff $A(c) := (L_c, R_c\iv, L_c R_c\iv)$ is an autotopism, and
$c$ is Moufang iff $B(c) := (L_c, R_c, R_c L_c)$ is an autotopism.
Since $c$ \emph{is} Moufang,
$D(c) := A(c) B(c) = (L_{c^2}, I, L_{c^2})$,
and $c$ is extra iff $D(c)$ is an autotopism iff $c^2 \in N$.

Now, (ii) holds iff $H(c) := (L_cR_c, L_c\iv, L_c)$ is an autotopism,
and $c$ Moufang implies that $F(c) := (R_c, L_c\iv R_c, R_c)$ is
an autotopism (see $({\Fex}2)$ in Lemma \lemref{eqns}). Then
(i)$\leftrightarrow$(ii) follows from $H(c) = A(c)F(c)$
(using $L_cR_c = R_c L_c$).  (i)$\leftrightarrow$(iii) is similar.
\end{proof}

\begin{lemma}
\lemlabel{ex-mfg}
In a CC-loop $Q$, let $a$ be an extra element and let $b$ be a Moufang element.
Then $ab$ and $ba$ are Moufang elements.
\end{lemma}

\begin{proof}
By Lemma \lemref{nuc-square}, $a$ is Moufang
and $a^2\in N(Q)$. By Lemma \lemref{Q/N-group}(1), 
$(x,a\riv,b) = (x, a^2 a\riv, b) = (x,a,b)$. Now use
Lemma \lemref{M-products}.
\end{proof}

\begin{lemma}
\lemlabel{square-subloop}
In a CC-loop $Q$, 
$S(Q) := \{ a : a^2 \in N(Q) \}$ is a normal subloop.
\end{lemma}

\begin{proof}
If $a,b\in Q$, then by Theorem \thmref{basarab},
there exists $n_1,n_2\in N$ such that
$(ab)^2 = (a^2 b^2) n_1$ and $(a\riv)^2 = (a^2)\riv n_2$.
Thus $a,b\in S$ implies $ab, a\riv \in S$, and so $S$
is a subloop by Lemma \lemref{subloop}.
Normality follows from Lemma \lemref{normal}.
\end{proof}

\begin{theorem}
\thmlabel{extra-sblp}
In a CC-loop $Q$, $Ex(Q)$ is a normal subloop.
\end{theorem}

\begin{proof}
Fix $a,b\in Ex(Q)$. By Lemmas \lemref{nuc-square} and \lemref{ex-mfg}, $ab$ is
Moufang.  By Lemmas \lemref{nuc-square} and \lemref{square-subloop},
$(ab)^2 \in N$, so $ab \in Ex(Q)$. Next
$a\iv \in Ex(Q)$ by Corollary \cororef{mfg-more}.
Thus $Ex(Q)$ is a normal subloop by
Lemmas \lemref{normal} and \lemref{subloop}.
\end{proof}

Now we relate Moufang and pseudoMoufang elements to WIP elements,
after first proving a technical lemma.

\begin{lemma}
\lemlabel{ps-char}
Let $c$ be an element of a CC-loop. The following are equivalent:
(i) $c$ is pseudoMoufang, (ii) $cx = (cx\cdot yx\riv)\cdot xy\riv$ for all $x,y$,
(iii) $xc = y\liv x\cdot (x\liv y\cdot xc)$ for all $x,y$.
\end{lemma}

\begin{proof}
$(ii)\to (i)$ follows from taking $x = 1$ and using Corollary \cororef{E-Mfg}(ii).

For $(i)\to (ii)$: We start with the equation $x = (x\cdot yx\riv) \cdot xy\riv$
(\cite{KKP}, Lemma 5.1). Define $u,v$ by the equations $yx\riv = (u)T_x$ and $xy\riv = (v)T_x$. Then 
$x = ux\cdot (v)T_x = uv\cdot x$ using \RCC. Thus $uv = 1$, and so by Corollary \cororef{E-Mfg}(ii),
$c = c E_u = cu\cdot v$. Hence $cx = (cu\cdot v)x = (cu\cdot x) \cdot xy\riv
= (cx \cdot yx\riv) \cdot xy\riv$ using \RCC\ twice.

The proof of $(i)\leftrightarrow (iii)$ is just the mirror of the preceding argument.
\end{proof}

\begin{lemma}
\lemlabel{wip-extra}
Let $c$ be an element of a CC-loop.
Any two of the following properties imply
the third: (i) $c$ is a WIP element, (ii) $c$ is Moufang,
(iii) $c$ is pseudoMoufang.
In case these conditions hold, $c$ is extra.
\end{lemma}

\begin{proof}
If $c$ is a WIP element, then (ii) and (iii) are
equivalent by Corollaries \cororef{wip2} and \cororef{E-Mfg}.
If $c$ is Moufang and pseudoMoufang, then
$c = c E_x = (c \cdot xE_c)\cdot x\riv$, and so $c$ is WIP by
Lemma \lemref{wip}.

Now suppose $c$ satisfies (i), (ii), and (iii). Since $c\iv$
is pseudoMoufang (Theorem \thmref{other-sblp}), we have
$c\iv x = (c\iv x\cdot yx\riv)\cdot xy\riv$ for all $x,y$,
by Lemma \lemref{ps-char}.  Replacing $x$ with $cx$ and using $c\iv \cdot cx = x$
(Corollary \cororef{E-Mfg}(i)), we obtain $x = [x\cdot y(cx)\riv ] (cx\cdot y\riv)$.
Now $(cx)\riv = x \ld c\iv$ by Lemma \lemref{wip}, and so 
$x = [x \cdot y (x \ld c\iv)] (cx\cdot y\riv )
= ((xy)\rd x) c\iv\cdot (cx\cdot y\riv )$, using \LCC.
Thus $xc = [((xy)\rd x) c\iv \cdot c] (c\ld [(cx\cdot y\riv )c])
= ((xy)\rd x) (c\ld [(cx\cdot y\riv )c])$ by \RCC\ and Corollary \cororef{E-Mfg}(1)
again. Using \LCC\ again,
$xc = x\cdot y \{ x\ld (c\ld [(cx\cdot y\riv )c])\}$. Cancelling $x$'s and rearranging, we have
$(cx\cdot y\riv)c = c (x(y\ld c)) = c(x\cdot y\riv c)$ by Corollary \cororef{E-Mfg}(2).
Replacing $y$ with $y\liv$ establishes the desired result.
\end{proof}

\begin{example}
\exlabel{extra-notwip}
The converse of Lemma \lemref{wip-extra} is not true, since in a
CC-loop, extra elements need not be WIP.  For the CC-loop $Q$
in Table \tableref{extra},
$Z(Q) = \{0,1\}$, $N(Q) = W(Q) = \{0,1,2,3\}$, 
$Ex(Q) = M(Q) = \{0,1,2,3,4,5,6,7\}$, and $P(Q) = \{0, 1, 2, 3, 8, 9, 10, 11\}$.
For example, $4$ is extra but not WIP because
$4 \cdot (8 \cdot4)\riv = 4 \cdot 14\riv = 4 \cdot 13 = 9 \ne 8\riv = 8$.
The set of power-associative elements of this loop is 
$\{0, 1, 2, 3, 4, 5, 6, 7, 8, 9, 10, 11\}$, which is not a subloop.
We also see that Theorem \thmref{cubes} cannot be improved;
that is, in a PACC-loop, $W(Q)$ contains all cubes, but in this loop,
$4$ is power-associative, and $4^3 = 4 \notin W(Q)$.
This example was produced by the program Mace4 \cite{McM};
as usual, given the example, it is trivial to verify its properties
using a standard programming language.
\end{example}
\begin{table}[htb]
{ 
\small
\setlength {\arraycolsep}{2pt}
\renewcommand {\arraystretch}{1.1}
\newcommand \sk{\ \;}
\newcommand\spx{\rule[-1pt]{0mm}{14pt}}  
\[
\begin{array}{c|cccc @{\sk} cccc @{\sk} cccc @{\sk} cccc}
    Q & 0& 1& 2& 3& 4& 5& 6& 7& 8& 9&10&11&12&13&14&15 \\
\hline
    0 & 0& 1& 2& 3& 4& 5& 6& 7& 8& 9&10&11&12&13&14&15 \\
    1 & 1& 0& 3& 2& 5& 4& 7& 6& 9& 8&11&10&13&12&15&14 \\
    2 & 2& 3& 0& 1& 6& 7& 4& 5&10&11& 8& 9&14&15&12&13 \\
    3 & 3& 2& 1& 0& 7& 6& 5& 4&11&10& 9& 8&15&14&13&12 \\
    4 & 4& 5& 6& 7& 0& 1& 2& 3&12&13&14&15& 8& 9&10&11 \spx \\
    5 & 5& 4& 7& 6& 1& 0& 3& 2&13&12&15&14& 9& 8&11&10 \\
    6 & 6& 7& 4& 5& 2& 3& 0& 1&14&15&12&13&10&11& 8& 9 \\
    7 & 7& 6& 5& 4& 3& 2& 1& 0&15&14&13&12&11&10& 9& 8 \\
    8 & 8& 9&11&10&14&15&13&12& 0& 1& 3& 2& 6& 7& 5& 4 \spx \\
    9 & 9& 8&10&11&15&14&12&13& 1& 0& 2& 3& 7& 6& 4& 5 \\
   10 &10&11& 9& 8&12&13&15&14& 2& 3& 1& 0& 4& 5& 7& 6 \\
   11 &11&10& 8& 9&13&12&14&15& 3& 2& 0& 1& 5& 4& 6& 7 \\
   12 &12&13&15&14&10&11& 9& 8& 5& 4& 6& 7& 3& 2& 0& 1 \spx \\
   13 &13&12&14&15&11&10& 8& 9& 4& 5& 7& 6& 2& 3& 1& 0 \\
   14 &14&15&13&12& 8& 9&11&10& 7& 6& 4& 5& 1& 0& 2& 3 \\
   15 &15&14&12&13& 9& 8&10&11& 6& 7& 5& 4& 0& 1& 3& 2
\end{array}
\]
}
\caption{extra $\not\rightarrow$ WIP}
\tablelabel{extra}
\end{table}

In the power-associative case, however, the converse
of Lemma \lemref{wip-extra} does hold, and the
pseudoMoufang elements coincide with the extra elements.

\begin{lemma}
\lemlabel{ex-pacc}
If $Q$ is a PACC-loop, then $P(Q) = Ex(Q)\leq W(Q)$.
\end{lemma}

\begin{proof}
For $Ex(Q) \leq W(Q)$: 
if $c\in Q$ is extra, then $c^2\in N(Q)$ (Lemma \lemref{nuc-square})
and $c^3\in W(Q)$ (Theorem \thmref{cubes}), and so
$c = c^3\cdot c^{-2} \in W(Q)$ by Theorem \thmref{W-subloop}.

For $P(Q) \leq W(Q)$:
if $c\in Q$ is pseudoMoufang, then
$cx\cdot x\iv = c = cx \cdot (cx)\iv c$
for all $x$ by Corollary \cororef{E-Mfg}.
Cancel the $cx$ to get that $c$ is WIP.

Finally, $Ex(Q) = P(Q)$ follows from Lemmas
\lemref{wip-extra} and \lemref{nuc-square}.
\end{proof}

\begin{corollary}
\corolabel{extra-group}
Let $a,b$ be elements of a PACC-loop, and suppose $a$
is extra. Then $\sbl{a,b}$ is a group.
\end{corollary}

\begin{proof}
By Lemmas \lemref{nuc-square}, \lemref{ex-pacc},  and \lemref{wip-extra},
$a$ is both Moufang and pseudoMoufang. Now apply Corollary \cororef{mfg-group}.
\end{proof}

We insert here some criteria for
determining when a PACC-loop is WIP.

\begin{lemma}
\lemlabel{wipequivs}
Let $Q$ be a PACC-loop.
The following are equivalent:
(i) $Q$ is WIP, (ii) every square is nuclear,
(iii) every square is extra, (iv) every square
is Moufang, (v) every square is pseudoMoufang,
(vi) every square is WIP.
\end{lemma}

\begin{proof}
(i)$\rightarrow$(ii) is due essentially to
{\cyr Basarab} \cite{BAS2},
and is true in any CC-loop; see also \cite{KKP}, \S7.
(ii)$\rightarrow$(iii) is clear from the definitions,
and (iii)$\rightarrow$(iv) follows from Lemma \lemref{nuc-square}.

(iv)$\rightarrow$(v):  By (iv) and Lemma \lemref{Ebasic}(4),
$E_y^4 = E_{y^2} = I$.  Since $\sbl{x^2, y^3}$ is a group
(by Corollary \cororef{old-groups}),
$x^2 E_y =  x^2 E_y^9 = x^2 E_{y^3} =  x^2 $, so $x^2$ is pseudoMoufang.

(v)$\rightarrow$(vi) is from Lemma \lemref{ex-pacc}.

(vi)$\rightarrow$(i): 
Each $x^3\in W(Q)$ by Theorem \thmref{cubes}. Thus
each $x = x^3\cdot x^{-2} \in W(Q)$ by Theorem \thmref{W-subloop}.
\end{proof}

We now turn to the main results of this section.

\begin{theorem}
\thmlabel{x4}
In a CC-loop $Q$, if $c\in W(Q)$, then $c^2\in Ex(Q)$ and
$c^2 c^2\in N(Q)$.
\end{theorem}

\begin{proof}
For each $x\in Q$, $c^2\in N(\sbl{c,x})$ by Lemma \lemref{wip-2gen}.
Thus $c^2$ is Moufang, and so by Lemma \lemref{wip-extra},
$c^2$ is extra. The rest follows from Lemma \lemref{nuc-square}.
\end{proof}

Incidentally, by Lemma \lemref{wip-2gen},
$c\cdot cc^2 = c^2c^2 = c^2c\cdot c$;
however, $c \cdot c^2$ need not equal $c^2 \cdot c$.
Thus, the final part
of the conclusion of Theorem \thmref{x4} would have been
ambiguous if written as $c^4\in N(Q)$.

\begin{theorem}
\thmlabel{x12}
For all $c$ in a PACC-loop $Q$,
\begin{enumerate}
\item[1.] \qquad $c^3$ is a WIP element
\item[2.] \qquad $c^6$ is an extra element
\item[3.] \qquad $c^{12} \in N(Q)$
\end{enumerate}
In particular, if $c$ has finite order prime to $6$, then $c\in N(Q)$.
\end{theorem}

\begin{proof}
(1) just restates Theorem \thmref{cubes}, and apply 
Theorem \thmref{x4} for the rest. 
\end{proof}

Note that Theorem \thmref{x12}(2) and Corollary \cororef{extra-group}
give a different proof of part of Corollary \cororef{old-groups}.
We now have:

\begin{corollary}
\corolabel{quotients}
Let $Q$ be a PACC-loop.
\begin{enumerate}
\item[1.] \qquad $Q/W(Q)$ is an elementary abelian $3$-group.
\item[2.] \qquad $Q/Ex(Q)$ is an abelian group of exponent $6$.
\item[3.] \qquad $Q/N(Q)$ is an abelian group of exponent $12$.
\end{enumerate}
\end{corollary}

\section{Fat Nuclei}
\seclabel{fat}

Following up on Corollary \cororef{quotients}(3), in this
section we consider the minimal possibilities for $Q/N(Q)$
for PACC-loops $Q$.

\begin{lemma}
\lemlabel{GN}
If $Q$ is a PACC-loop, $G \le Q$ is a group,
and $GN = Q$, then $Q$ is a group.
\end{lemma}
\begin{proof}
The set $G \cup N$ associates, so $\sbl{G \cup N}$ is a group;
see (\cite{KKP}, Corollary 6.4).
\end{proof}

\begin{corollary}
\corolabel{Q/N-cyclic}
A PACC-loop $Q$ with $Q/N$ cyclic is a group.
\end{corollary}

The following two lemmas are especially useful for PACC-loops
in which the center coincides with the nucleus.

\begin{lemma}
\lemlabel{nzcomp}
Let $Q = \sbl{a,b}$ be a PACC-loop, and assume $bE_a = bu$ and
$aE_b = av$ where $u,v\in Z(Q)$. Then
$u^3 = v^3$, $u^6 = v^6 = 1$, and 
\[
(a^i b^j, a^k b^{\ell}, a^m b^n) =
u^{-ikn-i{\ell}m-jkm} v^{-i{\ell}n-jkn-j{\ell}m}
\tag{$\dagger$}
\]
for all integers $i,j,k,\ell,m,n$. Also,
$a^6, b^6 \in N(Q)$ and $A(Q) = \sbl{u,v} \leq Z(Q)$.
\end{lemma}

\begin{proof}
Since $E_a^6 = E_b^6 = I$ (Lemma \lemref{Ebasic}(5)), 
we have $u^6 = v^6 = 1$. Now $u = (b,a,a\iv)
= (a,a,b)\iv$ and $v = (a,b,b\iv) = (a,b,b)\iv$
by Lemmas \lemref{Ebasic}(2) and \lemref{assoc-props}(3).
Applying Lemma \lemref{assoc-props}(1)(3) and multiple
inductions, we get
\begin{align*}
(a^i b^j, a^k b^{\ell}, a^m b^n) &=
(a,a,b)^{ikn+i{\ell}m+jkm} (a,b,b)^{i{\ell}n+jkn+j{\ell}m} \\
&= u^{-ikn-i{\ell}m-jkm} v^{-i{\ell}n-jkn-j{\ell}m}.
\end{align*}
Now if $i=j=k={\ell}=m=n=1$, we have
$1 = (ab,ab,ab) = u^{-3} v^{-3}$, and when
combined with $u^6 = v^6 = 1$, this gives $u^3 = v^3$.
By Theorem
\thmref{basarab} and Lemmas \lemref{T-comm} and
\lemref{Q/N-group}, every associator in $Q$ can be
expressed in the form $(\dagger)$, and so $A(Q) = \sbl{u,v}\leq Z(Q)$. 
Taking $i = 6$, $j=0$ in $(\dagger)$ gives $a^6\in N(Q)$,
and similarly, $b^6\in N(Q)$.
\end{proof}

\begin{corollary}
\corolabel{nz-coro}
Let $Q$ be a $2$-generated PACC-loop with $A(Q)\leq Z(Q)$.
Then $Q/N(Q)$ is an abelian group of exponent $6$.
\end{corollary}

\begin{lemma}
\lemlabel{nzcomp2}
Let $Q = \sbl{a,b}$ be a PACC-loop satisfying the hypotheses
of Lemma \lemref{nzcomp}.
In addition, assume $ba = abz$ where $z = [b,a]\in Z(Q)$. Then
\[
b^i a^j = a^j b^i \cdot z^{ij} u^{i(j-1)j/2} v^{-j(i-1)i/2}\ ,
\tag{Z1}
\]
and
\[
a^i b^j  \cdot a^k b^\ell =
 a^{i + k} b^{j+\ell} \cdot z^{jk} u^{i\ell k + j(k-1)k/2} v^{-ij\ell -k(j-1)j/2}
\tag{Z2}
\]
for all integers $i,j,k,\ell$.
\end{lemma}

\begin{proof}
We first prove the following special cases of (Z1):
\[
b a^j = a^j b \cdot z^j u^{(j-1)j/2} \qquad \mathrm{(Z3)} \qquad\qquad\qquad
b^i a = a b^i \cdot z^i v^{-(i-1)i/2}\qquad  \mathrm{(Z4)}
\]
Now (Z3) is clear for $j = 1$, and proceeding by induction,
using Lemma \lemref{nzcomp}:
\begin{align*}
b a^{j+1} &= ba\cdot a^j\cdot (b,a,a^j)\iv 
= abz\cdot a^j \cdot u^j = a\cdot ba^j \cdot (a,b,a^j)\cdot zu^j \\
&= a\cdot a^j b \cdot z^j u^{(j-1)j/2} \cdot (a,a^j,b)\cdot zu^j
= a^{j+1} b \cdot z^{j+1} u^{j(j+1)/2}
\end{align*}
The mirror of this argument yields (Z4).
Now, in (Z3), we can replace $b$ by $b^i$,
replace $z$ by $z^i v^{-(i-1)i/2}$ (using (Z4)),
and replace $u$ by $u^i$ (since $(b^i) E_a = b^i u^i$)
to get (Z1).

For (Z2), we repeatedly use Lemma \lemref{nzcomp}, (Z1) and (Z2)
to compute 
\begin{align*}
a^i b^j  \cdot a^k b^{\ell} &=
a^i \cdot ( b^j  \cdot a^k b^{\ell} ) \cdot (a^i, b^j, a^k b^{\ell}) \\
&= a^i (b^j a^k\cdot b^{\ell}) \cdot (b^j, a^k, b^{\ell})\iv \cdot u^{-ijk} v^{-ij\ell} \\
&= a^i (a^k b^j\cdot b^{\ell}) \cdot
z^{jk} u^{j(k-1)k/2} v^{-k(j-1)j/2} \cdot
(b^j, a^k, b^{\ell})\iv \cdot
u^{-ijk} v^{-ij\ell} \\
&= a^i \cdot a^k b^{j+\ell}\cdot 
z^{jk} u^{-ijk+j(k-1)k/2} v^{-ij\ell -k(j-1)j/2} \\
&= a^{i+k} b^{j+\ell}\cdot (a^i,a^k,b^{j+\ell})\iv \cdot
z^{jk} u^{-ijk+j(k-1)k/2} v^{-ij\ell -k(j-1)j/2} \\
&= a^{i + k} b^{j+\ell} \cdot z^{jk} u^{i\ell k+j(k-1)k/2} v^{-ij\ell -k(j-1)j/2}
\end{align*}
\end{proof}

In Theorem \thmref{f-good}, we will construct loops satisfying
the hypotheses of Lemma \lemref{nzcomp2}.

\begin{lemma}
\lemlabel{nz3}
Let $Q = \sbl{a,b}$ be a PACC-loop satisfying the hypotheses of
Lemma \lemref{nzcomp2}, and assume also that $a^3\in Z(Q)$.
Then $u^3 = v^3 = z^3 = 1$ and $b^3\in Z(Q)$.
\end{lemma}

\begin{proof}
By Lemma \lemref{nzcomp2}(Z1), if $a^3\in Z(Q)$, then 
$1 = z^{3i} u^{3i} v^{-3(i-1)i/2}$ for all $i$. Taking
$i = 1$, we have $1 = z^3 u^3$. Since $u^6 = 1$,
$z^3 = u^3$. Then taking $i = 2$, we
get $1 = z^6 u^6 v^{-3} = v^{-3}$, and so
$v^3 = 1$, and so $u^3 = z^3 = 1$ by Lemma \lemref{nzcomp}. That $b^3\in Z(Q)$
then follows from Lemma \lemref{nzcomp2}(Z2).
\end{proof}

\begin{lemma}
\lemlabel{nz2}
Let $Q = \sbl{a,b}$ be a PACC-loop satisfying the hypotheses of
Lemma \lemref{nzcomp2}, and assume also that $a^2, b^2 \in Z(Q)$.
Then $u = v = z^2$ and $u^2 = z^4 = 1$.
\end{lemma}

\begin{proof}
By Lemma \lemref{nzcomp2}(Z1), if $a^2\in Z(Q)$, then 
$1 = z^{2i} u^i v^{-(i-1)i}$ for all $i$. Taking
$i = 1$, we have $1 = z^2 u$, and so $u = z^{-2}$.
Taking $i=2$, we get $1 = z^4 u^2 v^{-2} = v^{-2}$.
Next, $b^2\in Z(Q)$ gives $1 = z^{2j} u^{(j-1)j} v^{-j}$ for all $j$.
Taking $j = 1$ gives $1 = z^2 v\iv = z^2 v$, and so
$u = z^{-2} = v$ and $u^2 = z^{-4} = v^2 = 1$.
\end{proof}

Since $[a,b] \in N(Q)$ always holds in a CC-loop
(Theorem \thmref{basarab}), we have:

\begin{corollary}
\corolabel{NZ-bool}
Let $Q$ be a PACC-loop.
\begin{enumerate}
\item[1.] If $Z(Q)$ is an elementary abelian $2$-group
and $A(Q) \leq Z$, then $Q$ is WIP. 
\item[2.] If $N = Z$ is an elementary abelian $2$-group, then $Q$ is extra.
\end{enumerate}
\end{corollary}

\begin{proof}
For (1): 
fix $a,b\in Q$ and adopt the notation of Lemma \lemref{nzcomp}.
Then $u = v$ and $u^2 = v^2 = 1$. By Lemma \lemref{Ebasic}(4), $b E_{a^2} = b E_a^4 = b$,
and so every square is Moufang.
By Lemma \lemref{wipequivs}, $Q$ is WIP.

For (2): $N = Z$ being an elementary abelian $2$-group, and
Theorem \thmref{basarab} imply
the hypotheses of (1), so that $Q$ is WIP. Thus squares
are central by Lemma \lemref{wipequivs}. Fixing $a,b\in Q$,
Lemma \lemref{nz2} implies $u = v = z^2 = 1$. Hence $\{a,b\}$
associates so that $\sbl{a,b}$ is a group, \emph{i.e.}, $Q$ is extra.
\end{proof}

\begin{corollary}
\corolabel{order8}
If $Q$ is a PACC-loop of order $8$, then $Q$ is a group.
\end{corollary}
\begin{proof}
Otherwise, by Corollaries \cororef{Q/N-cyclic} and \cororef{NZ-bool},
$|N| = 2$ and $Q$ is an extra loop, but the
smallest nonassociative extra loop has order $16$ (see \cite{KK, KKP},
or \cite{CH}).
\end{proof}

We conclude this section by examining the case
$|Q/N| = 4$ in some detail.

\begin{lemma}
\lemlabel{index4}
Assume that $Q$ is a PACC-loop with $|Q/N| = 4$.
Then $Q/N$ is an elementary abelian $2$-group, and $Q$ has WIP.
If $Q/N = \{N, Na, Nb, Nab\}$, then
$A(Q) = \{1,u\} \leq Z(Q)$, where $u = (a,a,b) = (a,b,b) \ne 1$.
Also, $4$ divides $|N|$.
\end{lemma}
\begin{proof}
$Q/N$ is an elementary abelian $2$-group by Corollary \cororef{quotients}(1)
(or by Lemma \lemref{wip-2gen}),
and $Q$ has WIP by Corollary \cororef{nucwip}.
The claim about $A(Q)$ follows from Lemma \lemref{wip-2gen}.

Now suppose that $|N| = 2r$, where $r$ is odd,
so $|Q| = 8r$.  Let $G = Q/A$, and let $\pi : Q \twoheadrightarrow G$
be the quotient map.  $G$ is a group of order $4r$
and $\pi(N) \normal G$, with $|\pi(N)| = r$.
Let $P$ be a Sylow 2-subgroup of $G$.  Then $|P| = 4$,
and $P$ contains exactly one element from each of the four cosets of $\pi(N)$.
Say $P = \{1, \pi(n_1 a), \pi(n_2 b), \pi(n_3 ab)\}$, with
$n_1, n_2, n_3 \in N$.  Then $\pi\iv(P)$ is a subloop of $Q$ of order $8$.
Since $(n_1 a,n_1 a,n_2 b) = (a,a,b) = u \ne 1$,
$\pi\iv(P)$ is nonassociative, contradicting Corollary \cororef{order8}.
\end{proof}

\section{Orders}
\seclabel{orders}
There are nonassociative PACC-loops of order $16$
(the five extra loops plus three others;
see \S\secref{2loops}) and of order $27$ (see
\S\secref{27}).  By taking products with a group 
of order $n$, we get nonassociative PACC-loops of order $16 n$
and $27 n$ for all finite $n \ge 1$.  These are the only possible
finite orders by Theorem \thmref{16or27} below.

\begin{lemma}
\lemlabel{N-prime}
Let $Q$ be a PACC-loop with finite associator subloop $A = A(Q)$.
Let $n$ be an integer relatively prime to $|A|$.
\begin{enumerate}
\item[1.] If $a\in Q$ satisfies $E_a^n = I$, then $E_a = I$.
\item[2.] If $a\in Q$ satisfies $E_{a^n} = I$, then $E_a = I$.
\end{enumerate}
\end{lemma}

\begin{proof}
For (1): For each $x\in Q$, $(x)E_a = xu$ where $u = (x,a,a\iv)$.
Thus $x = (x)E_a^n = xu^n$, and $u^n = 1$. Since $n$
is prime to $|A|$, $u = 1$.

For (2): this follows from (1) and $E_{a^n} = E_a^{n^2}$.
\end{proof}

\begin{lemma}
\lemlabel{NQ-prime}
Let $Q$ be a finite PACC-loop, and suppose $|A(Q)|$ is
relatively prime to $|Q/N(Q)|$. Then $Q$ is a group.
\end{lemma}

\begin{proof}
For each $a \in Q$, we have $a^k \in N$
for some $k$ relatively prime to $|A|$.
Thus $E_{a^k} = I$, and so $E_a = I$
by Lemma \lemref{N-prime}. Therefore 
$Q$ is Moufang (Corollary \cororef{E-Mfg}(i)), and hence 
an extra loop. But in an extra loop, $A$ and
$Q/N$ are elementary abelian $2$-groups \cite{KK}, and so
$Q$ must be a group.
\end{proof}

\begin{lemma}
\lemlabel{WIP3N}
If  $Q$ is a finite 
PACC-loop and $3 \nmid |A(Q)|$, then $Q$ is WIP.
\end{lemma}
\begin{proof}
First, fix any $x,y \in Q$.  Then $y E_x = y n$ for some $n \in A$.
$E_x^6 = I$ implies $n^6 = 1$, but then $n^2 = 1$ since
 $3 \nmid |A|$.
Thus, $E_x^2 = I$ for all $x$, and hence also $E_{x^2} = E_x^4 = I$.
Thus, $Q$ is WIP by Lemma \lemref{wipequivs} and Corollary \cororef{E-Mfg}.
\end{proof}

\begin{theorem}
\thmlabel{16or27}
Let $Q$ be a finite, nonassociative
PACC-loop.
If $Q$ is WIP, then $16 \mid |Q|$.
If $Q$ is not WIP, then $27 \mid |Q|$.
\end{theorem}

\begin{proof}
If $Q$ is WIP, then $x^4 \in N$ for all $x$ by Theorem \thmref{x4}.
Thus, $|Q/N| = 2^k$ for some $k \ge 2$ by Corollary \cororef{Q/N-cyclic},
and $2 \mid |N|$ by Lemma \lemref{NQ-prime},
so $16 \mid |Q|$ unless $4 \nmid |N|$ and $|Q/N| = 4$;
but this would contradict Lemma \lemref{index4}.

If $Q$ is not WIP, fix $b \notin W$.
Say  $o(b) = 3^j k$, where $\gcd(3,k) = 1$.
Then $b^k \notin W$ by Corollary
\cororef{Q/W}, and so replacing $b$ by 
$b^k$ if necessary, WLOG assume that $o(b) = 3^j$.
Now fix $a$ such that $b (ab)\iv \ne a\iv$.
Then $a \notin \sbl{b} W$.
To see this, suppose that $a = b^i w $ for $w \in W$ and $i \in \ZZ$.
Then $a,b \in \sbl{w, b^2}$
(since $o(b) = 3^j$), and $\sbl{w, b^2}$ is a group
by Lemma \lemref{wip-2gen}, contradicting $b (ab)\iv \ne a\iv$.

Thus, $9 \mid |Q/W|$.  Since $3 \mid |A|$ by Lemma
\lemref{WIP3N}, $3\mid |W|$, and so $27 \mid |Q|$.
\end{proof}

\section{Extension}
\seclabel{ext}
Here, we show how in some cases, the equations of
Lemma \lemref{nzcomp} may be used as a prescription for constructing
a PACC-loop.  This will be useful in Sections
\secref{2loops} and \secref{27},
where we construct all PACC-loops of a given order.
The natural converse to Lemma \lemref{nzcomp2} would be:

\begin{lemma}
\lemlabel{nzcomp-rev}
Assume that $Q$ is a PACC-loop, $a,b \in Q$, and $z,u,v \in Z(Q)$, with
\[
a^i b^j  \cdot a^k b^\ell =
 a^{i + k} b^{j+\ell} \cdot z^{jk} u^{j(k-1)k/2} v^{-k(j-1)j/2}
\cdot u^{i \ell k} v^{-ij \ell}  
\]
holding for all $i,j,k,\ell \in \ZZ$.  Then
$ba = ab z$, $(b)E_a = b u$, and
$(a) E_b = a v$.
\end{lemma}
\begin{proof}
Setting $j = k = 1$ and $i = \ell = 0$, we get $ba = ab z$.
Also,
\[
(b)E_a = b a \cdot a\iv = ab \cdot a\iv \cdot z   =
b \cdot z^{-1} u  \cdot z = bu 
\]
(setting $i = j = 1, k = -1, \ell = 0$),
and
\[
(a)E_b = a b \cdot b\iv =  a \cdot  v
\]
(setting $i = j = 1, k = 0, \ell = -1$).
\end{proof}

We now consider how to build such loops.
First, a more general construction; the following is a variant
of the construction described in \cite{KK} (see Definition 7.1
and Lemmas 7.2 and 7.3):

\begin{definition}
\deflabel{extn}
Say we are given abelian groups $(A,+)$ and $(G,+)$, and 
a function $f : A \times A \to G$, satisfying
\[
f(0,a) = f(a,0) = 0 \quad \text{for all}\  a\in A \ \ .
\eqno{(*1)}
\]
Then $A \ltimes_f G$ is the set $Q = A \times G$ with a product 
defined by:
\[
(a, x) \cdot (b, y) = (a + b ,\  x + y + f(a,b))  \ \ .
\]
\end{definition}

\begin{lemma}
\lemlabel{extn-basic}
Let $Q = A \ltimes_f G$, where $f$ satisfies $(*1)$.  Then
$Q$ is a loop with identity element $(0,0)$ and divisions given by:
\begin{align*}
(c,z) \rd (b,y) &= (c - b , \  z - y - f(c-b,b))  \\
(a,x) \ld (c,z) &= ( c - a , \ z - x - f(a,c-a) ) \ \ .
\end{align*}
Associators are given by
\[
( (a,x) , (b,y) , (c,z) ) = \big(0,\ \AAA_f(a,b,c) \big) \ \ ,
\]
where
\[
\AAA_f(a,b,c) := f(a,b) + f(a+b, c) - f(b,c) - f(a, b+c) .
\]
Commutators are given by
\[
\lbrack (b,y),(a,x) \rbrack = (0, f(b,a) - f(a,b))\ \ ,
\] 
and so
\[
(b,y) T_{(a,x)} = (b,y) \cdot (0, f(b,a) - f(a,b)) \ \ .
\]
Also, $G \cong \{0\} \times G \le Z(Q)$.
\end{lemma}

\begin{proof}
The division formulas are obtained by solving
$ (a, x) \cdot (b, y) =  (c,z)  $
for $(a,x)$ and for $(b,y)$.
It is clear that elements of $\{0\} \times G $ commute with all elements
of $Q$.  To show that $\{0\} \times G \le N(Q)$ (and hence
$\{0\} \times G \le Z(Q)$),  note
\begin{align*}
(a,x) \cdot (b,y) (c,z) &= (a + b + c,\  x + y + z + f(b,c) + f(a, b+c))\\
(a,x) (b,y) \cdot (c,z) &= (a + b + c,\  x + y + z + f(a,b) + f(a+b, c))\ \ .
\end{align*}
These two are equal if at least one of $a,b,c$ is $0$.  Thus,
$\{0\} \times G \le N(Q)$. 
Dividing the expressions for products, we get the formula
for associators.  Finally,
\[
(b,y) \cdot (a,x) = (a,x) \cdot (b,y) \cdot (0, f(b,a) - f(a,b)) \ \ .
\]
yields the formulas for commutators and for the mappings
$T_{(a,x)} = R_{(a,x)} L_{(a,x)}\iv$.
\end{proof}

\begin{definition}
\deflabel{good-f}
If $A,G$ be abelian groups, then a mapping $f : A\times A\to G$ is
\begin{itemize}
\item \emph{CC-good} if $f$ satisfies $(*1)$ and $A\ltimes_f G$ is a CC-loop.
\item \emph{PACC-good} if $f$ satisfies $(*1)$ and $A\ltimes_f G$ is a PACC-loop.
\end{itemize}
\end{definition}

The characterization of associators in Lemma \lemref{extn-basic}
gives us:

\begin{lemma}
\lemlabel{extn-cc}
Let $A,G$ be abelian groups, and assume $f: A\times A\to G$ satisfies
$(*1)$. Then $f$ is CC-good iff $\AAA_f(a,b,c)$ is invariant under
permutations of $\{a,b,c\}$.
In that case, $(a,x) \in Z(A\ltimes_f G)$ iff $f(a,b) = f(b,a)$ for all $b$.
\end{lemma}

\begin{proof}
Since $A \ltimes_f G$ clearly satisfies (b) of Lemma \lemref{ccequiv},
CC is equivalent to (a), which is equivalent to the
invariance of  $\AAA_f(a,b,c)$ under permutations.
Then, note that for CC-loops, $(a,x) \in Z(Q)$ iff $(a,x)$ commutes
with all other elements, and apply the characterization of commutators
in Lemma \lemref{extn-basic}.
\end{proof}

\begin{lemma}
\lemlabel{extn-pacc}
Let $A,G$ be abelian groups, and assume $f: A\times A\to G$ is CC-good.
Then $f$ is PACC-good iff $\AAA_f(a,a,a) = 0$ for all $a \in A$.
In that case,
$f(ma,na) =   f(na,ma)$ for all $a \in A$ and all $m,n \in \ZZ$.
\end{lemma}

\begin{proof}
By Lemma \lemref{PA-char}, $(x,a)$ is a power-associative element of $A\ltimes_f G$ iff
$(0,0) = ((x,a),(x,a),(x,a)) = (0,\AAA_f(a,a,a))$.
The rest follows since $(a,x)^m$ and $(a,x)^n$ must commute
in a power-associative loop.
\end{proof}

\begin{lemma}
\lemlabel{good-add}
If $A,G$ are abelian groups, then the sets of all CC-good and
PACC-good $f:A\times A\to G$ are subgroups of $G^{A \times A}$.
\end{lemma}

\begin{proof}
If $f,g:A\times A\to G$ satisfy $(*1)$, then obviously so does $f+g$.
The rest follows from Lemmas \lemref{extn-cc}, \lemref{extn-pacc},
and the observation that $\AAA_{f+g}(a,b,c) = \AAA_f(a,b,c) + \AAA_g(a,b,c)$
for all $a,b,c\in A$.
\end{proof}

\begin{lemma}
\lemlabel{bilin}
If $A,G$ are abelian groups and $f: A \times A \to G$ is 
bilinear, then $f$ is PACC-good and $A \ltimes_f G$ is a group.
\end{lemma}

Motivated now by Lemma \lemref{nzcomp-rev}, we start with
$A = \sbl{a} \times \sbl{b} \cong \ZZ\times\ZZ$,
and let $G$ be some abelian group containing elements $z,u,v$.
Converting to additive notation, we define 
$f : A \times A \to G$ by:
\[
f(ia + jb , ka + \ell b) = (j(k-1)k/2 + i \ell k) u
+ (-k(j-1)j/2 -ij \ell) v  + (jk) z\ \ .
\eqno{(*2)}
\]
This is a well-defined function since $a$ and $b$ are assumed to have
infinite order. Then, if it is desired to construct a finite loop, we
will quotient out a suitable normal subloop of $A\ltimes_f G$.

\begin{theorem}
\thmlabel{f-good}
Let $G$ be an abelian group and
$A = \sbl{a} \times \sbl{b} \cong \ZZ\times\ZZ$.
Fix $z,u,v\in G$ and define $f$ as in $(*2)$,
and let $P = A\ltimes_f G$.
Then $P$ is a CC-loop.
In addition, $P$ is a PACC-loop iff $3u = 3v$ and $6u = 6v = 0$.
In this case, $6a, 6b\in N(P)$.
\end{theorem}

\begin{proof}
First, observe that the mapping $f_z(ia+jb,ka+\ell b) = (jk)z$
is bilinear, and hence PACC-good by Lemma \lemref{bilin}.

Next, consider $f_u(ia+jb,ka+\ell b) = (j(k-1)k/2 + i \ell k) u$.
We compute $\AAA_{f_u}$:
\begin{align*}
& \AAA_{f_u}(ia+jb,ka+\ell b,pa+qb)  = \\
& f_u(ia+jb,ka+\ell b) + f_u((i+k)a+(j + \ell) b, pa+qb) -  \\
& f_u(ka+\ell b,pa+qb) - f_u(ia+jb, (k + p)a+(\ell+q) b)  =  \\
& ( j(k-1)k/2 + i k \ell )u + 
 ( ( j+ \ell) p(p-1)/2 + (i + k) p q  )u - \\
& ( \ell p(p-1)/2 + k p q  )u - 
 ( j(k+p-1)(k+p)/2 + i (k+p) (\ell + q) )u  = \\
& (-ikq - jkp  - i \ell p)u \ \ .
\end{align*}
Thus $\AAA_{f_u}$ is invariant under permutation of its arguments, and
so $f_u$ is CC-good by Lemma \lemref{extn-cc}.

Similarly, for $f_v(ia+jb,ka+\ell b) = (-k(j-1)j/2 -ij \ell) v$,
we find that
$\AAA_{f_v}(ia+jb,ka+\ell b,pa+qb) = -i\ell p - j\ell q - jkq$,
and so $f_v$ is CC-good by Lemma \lemref{extn-cc}.

Since $f = f_u + f_v + f_z$, $f$ is CC-good by Lemma \lemref{good-add}.

Now $\AAA_f = \AAA_{f_u} + \AAA_{f_v}$, and so
$\AAA_f(ia+jb, ia+jb, ia+ jb) = -3 i^2 j u - 3 ij^2 v  = -3ij(iu + jv)$
for all $i,j$. If $f$ is PACC-good, then
$-3ij(iu + jv) = 0$ for all $i,j$ by Lemma \lemref{extn-pacc}.
Taking $i = -j = 1$, we have
$3u = 3v$, and so $3ij(i+j)u = 0$. Taking $i = j = 1$, we have
$6u = 6v = 0$. Conversely, if $6u = 0$, then since $ij(i+j)$ is
always even, $3ij(i+j)u = 0$, and so if $3u = 3v$, then
$3ij(iu + jv)  = 0$ for all $i,j$. Therefore $f$ is PACC-good
by Lemma \lemref{extn-pacc}.

Finally, note that 
$\AAA_f(6a,ka+\ell b,pa+qb) = -6(kq + lp)u -6(lq)v = 0$, and
so $6a \in N(P)$, and similarly $6b\in N(P)$.
\end{proof}

We now consider special cases of this construction.

\begin{corollary}
\corolabel{G-exp3}
Let $G$ be an abelian group and
$A = \sbl{a} \times \sbl{b} \cong \ZZ\times\ZZ$.
Define $f$ as in $(*2)$,
and let $P = A\ltimes_f G$. Then TFAE:
\begin{itemize}
\item[1.] For \emph{every} choice of $z,u,v \in G$, 
$P$ is a PACC-loop.
\item[2.] $G$ is of exponent $3$.
\end{itemize}
\end{corollary}

\begin{proof}
For (1)$\implies$(2), applying the theorem with $v = 2u$ gives
$3v = 0$ for all $v\in G$. The converse is clear.
\end{proof}

The other special case of Theorem \thmref{f-good} we will
consider is motivated by Lemma \lemref{nz2}.

\begin{corollary}
\corolabel{goodf2}
Let $G$ be an abelian group and
$A = \sbl{a} \times \sbl{b} \cong \ZZ\times\ZZ$.
Fix $z\in G$, set $v = 2z$, $u = -2z$, and
define $f$ as in $(*2)$. Let
$P = A\ltimes_f G$. Then $P$ is a PACC-loop iff
the order of $z$ divides $12$.
\end{corollary}

\begin{proof}
In this case, $\AAA_f(ia+jb, ia+jb, ia+ jb) = -3ij(iu + jv)
= 6ij(i-j)z$. Since $ij(i-j)$ is always even, $P$ is
power-associative iff $12z = 0$ by Lemma \lemref{extn-pacc}.
\end{proof}

We now apply the corollaries to get examples with
$a,b$ of finite order:

\begin{lemma}
\lemlabel{goodfquot}
Let $(G, \cdot)$ be an abelian group of exponent $3$,
and fix $z,u,v,t,w \in G$.
Then there is a PACC-loop $Q = \sbl{G \cup \{a,b\}}$ with 
$G \le Z(Q)$ and $Q/G \cong \ZZ_3 \times \ZZ_3$, such that $Q$ satisfies:
\[
ba = ab z \qquad (b)E_a = b u \qquad (a) E_b = a v \qquad
a^3 = t \qquad b^3 = w \ \ .
\]
\end{lemma}

\begin{proof}
Start with the loop $P = (\ZZ\times\ZZ)\ltimes_f G$
constructed in Theorem \thmref{f-good},
where $a,b$ have infinite order.
Then $(b,0)(a,0) = (a,0)(b,0) (0,z)$,  $(b,0)E_{(a,0)} = (b,0) (0,u) $,
and $(a,0) E_{(b,0)} = (a,0) (0,v)$ by Lemma \lemref{nzcomp-rev}.
Since $G$ has exponent 3, $(*2)$ yields
\[
f(ia + jb , ka + \ell b) =
(-jk^2 + i \ell k + jk) u  + (kj^2  -ij \ell - jk) v + (jk) z \ \ .
\eqno{(*3)}
\]
Also note:
\[
  f(3ia + 3jb , ka + \ell b) = f(ia + jb , 3ka + 3\ell b) = 0  \ \ ,
\]
so that all $(3ia + 3jb, x) \in Z(P)$ by Lemma \lemref{extn-cc}.
Let $H = \{(3ia + 3jb, -it -jw) : i,j \in \ZZ\}$.
Then $H \le Z(P)$, so $H$ is a normal subloop, and the lemma
is satisfied by $P/H$.
To see this, note that $f(ia,ka) = f(jb,\ell b) = 0$, so that
$(a,0)^3 = (3a, 0) \equiv (0,t) \pmod H$ and
$(b,0)^3 = (3b, 0) \equiv (0,w) \pmod H$.
\end{proof}

When $|G| = 3$, this lemma lists all PACC-loops of order 27;
see \S\secref{27}.

\begin{lemma}
\lemlabel{goodf2quot}
Let $(G, \cdot)$ be an abelian group,
and fix $z, t,w \in G$ with $z^4 = 1$.
Then there is a PACC-loop $Q = \sbl{G \cup \{a,b\}}$ with 
$G \le Z(Q)$ and $Q/G \cong \ZZ_2 \times \ZZ_2$, such that $Q$ satisfies:
\[
ba = ab z \qquad (b)E_a = b z^{-2} \qquad (a) E_b = a z^2 \qquad
a^2 = t \qquad b^2 = w \ \ .
\]
\end{lemma}

\begin{proof}
Start with the loop $P = (\ZZ\times\ZZ)\ltimes_f G$
constructed in Theorem \thmref{f-good},
where $a,b$ have infinite order.
Then $(b,0)(a,0) = (a,0)(b,0) (0,z)$,  $(b,0)E_{(a,0)} = (b,0) (0,-2z) $,
and $(a,0) E_{(b,0)} = (a,0) (0,2z)$ by Lemma \lemref{nzcomp-rev}.
Since $4z = 0$, $(*2)$ yields
\[
f(ia + jb , ka + \ell b) =
(2ik\ell + 2ij\ell -jk^2 -kj^2  -jk )z \ \ .
\]
Also note:
\[
 f(2ia + 2jb , ka + \ell b) = f(ia + jb , 2ka + 2\ell b) = 0   \ \ ,
\]
since $jk^2 + jk$ and $kj^2 + jk$ are always even;
thus all $(2ia + 2jb, x) \in Z(P)$ by Lemma \lemref{extn-cc}.
Let $H = \{(2ia + 2jb, -it -jw) : i,j \in \ZZ\}$.
Then $H \le Z(P)$, so $H$ is a normal subloop, and the lemma
is satisfied by $P/H$.
\end{proof}

When $G = \ZZ_4$, this lemma describes the two PACC-loops $Q$ of order
16 with $Z(Q) = N(Q) \cong \ZZ_4$;
see \S\secref{2loops}.

We conclude this section by considering the automorphic inverse
property.

\begin{definition}
\deflabel{aip}
A PACC-loop has the \textit{automorphic inverse property (AIP)} iff
the map $x \to x\iv$ is an automorphism.
\end{definition}

When $A \ltimes_f G$ is power-associative, $f(a,-a) = f(-a,a)$ for
all $a$ (see Lemma \lemref{extn-pacc}), and 
$(a,x)\iv = (-a, -x - f(a,-a))$.  Then,
\begin{align*}
( (a,x)\cdot (b, y)) \iv &=                       
   (-a-b, -x-y - f(a,b) - f(a+b, -a -b) )                      \\
(a,x) \iv \cdot (b, y)\iv &=                    
  (-a-b, -x-y + f(-a,-b) - f(a, -a) - f(b, -b)  )  \ \ ,
\end{align*}
so $A \ltimes_f G$ satisfies the AIP iff
\[
 f(a+b, -a -b) = f(a, -a) + f(b, -b) -  f(-a,-b)   - f(a,b)
\]
holds for all $a,b \in A$.  We remark that when $f$ is bilinear,
this reduces to $f(b,a) = f(a,b)$ which is equivalent to
the group  $A \ltimes_f G$ being abelian.
Replacing $a$ by $ia + jb$ and $b$ by $ka + \ell b$, we get
\begin{align*}
f(ia + jb +ka + \ell b, -(ia + jb  +ka + \ell b)) = \\
f(ia + jb , -(ia + jb )) + f(ka + \ell b, -(ka + \ell b))\\
-  f(-(ia + jb ),-(ka + \ell b))   - f(ia + jb ,ka + \ell b)
\end{align*}
Now, consider the case where $G$ has exponent $3$ and
$f$ is as in $(*3)$.  Then
\[
f(ia + jb , -(ia + jb)) =
(-ji^2 + i j i - ji) u  + (-ij^2  +ij j + ji) v - (ji) z =
-ij(u - v + z) \ \ .
\]
So, the requirement becomes
\begin{align*}
& -(i +k) (j + \ell) (u - v + z) = \\
& -ij(u + v + z) -k\ell (u - v + z) \\
& - [ ( jk^2 - i \ell k + jk) u  + (-kj^2  +ij \ell - jk) v + (jk) z    ]  \\
& - [ (-jk^2 + i \ell k + jk) u  + (kj^2  -ij \ell - jk) v + (jk) z    ] \ \ .
\end{align*}
This simplifies to
\[
(i \ell + k j) (u - v + z) =  2jk (u - v + z) 
\]
for all $i,j,k,\ell$, which is equivalent to $v = u+z$.
Passing to the quotient, we get

\begin{lemma}
\lemlabel{aip27}
If $Q$ is the loop constructed in Lemma \lemref{goodfquot},
then $Q$ has the AIP iff  $v = u+z$.
\end{lemma}
\begin{proof}
We had $Q = P/H$, where $P = (\ZZ\times\ZZ)\ltimes_f G$.
We have just seen that $v = u+z$ implies the AIP of $P$, and
hence of $P/H$.
Conversely, if $v \ne u+z$, then the AIP fails in $P$, with
$ [(ia + jb, 0) \cdot (ka + \ell b, 0)] \iv \,\ld\,
[(ia + jb, 0)\iv \cdot  (ka + \ell b, 0)\iv] = (0,z) \ne (0,0) $
for some $i,j,k,\ell$.  Since $(0,z) \notin H$, the 
AIP fails in $Q$ as well.
\end{proof}

This lemma will be used in \S\secref{27}
to determine which of the
PACC-loops of order 27 satisfy the AIP.

\section{2-Loops}
\seclabel{2loops}

In this section we show that there are eight
nonassociative PACC-loops of order $16$; this includes the five
extra loops already described by Chein (\cite{CH}, p.~49).
If $Q$ is such a loop, then $Z = Z(Q)$
is nontrivial, so our strategy is to analyze the various possibilities
for $Z$ and $N$.

\begin{lemma}
\lemlabel{basic2}
If $Q$ is a PACC-loop of finite order $2^n$, then
\begin{itemize}
\item[1.] \qquad $|Z(Q)| = 2^r$, where $0< r \le n$.
\item[2.] \qquad Q has WIP.
\item[3.] \qquad $Q / N(Q)$ is an elementary abelian $2$-group.
\end{itemize}
\end{lemma}
\begin{proof}
(1) is from \cite{KKP}, Cor.~3.5, and is true of all CC-loops.
(2) is by Corollary \cororef{Q/W}.
(3) is by Lemma \lemref{wipequivs}.
\end{proof}

The five nonassociative extra loops of order $16$ 
all have $Z(Q) = N(Q) \cong \ZZ_2$ and $Q/Z$ an elementary abelian $2$-group.
For the nonextra ones, we have two cases, described by:

\begin{lemma}
\lemlabel{16cases}
If $Q$ is a nonextra PACC-loop of order $16$, then $|N| = 4$,
$Q/N \cong \ZZ_2 \times \ZZ_2$, and either:
\begin{itemize}
\item[1.] \qquad $|Z| = 4$  and $Z \cong \ZZ_4$, or
\item[2.] \qquad $|Z| = 2$ and $|N| = 4$.
\end{itemize}
\end{lemma}

\begin{proof}
$|N| = 2$ or $|N| = 4$ by Corollary \cororef{Q/N-cyclic},
but $|N| = 2$ would contradict Corollary \cororef{NZ-bool}, 
which also yields $|Z| = 4 \to Z \cong \ZZ_4$.
\end{proof}

In Case (1), fix $a,b\in Q$ with
$aZ \neq bZ$. Then $Q / Z = \{ Z, aZ, bZ, abZ\}$.
Define $z = [b,a] = (ab) \ld (ba)$ as in Lemmas \lemref{nzcomp2} and \lemref{nz2}.
Then $z^2 \ne 1$, since otherwise 
$\{ a,b\}$ would associate, but then $Q$ would be a group by 
Lemma \lemref{GN}.  So, $z$ is a generator of $Z$. 
Say  $a^2 = z^r$ and $b^2 = z^s$, where $r,s \in \{0,1,2,3\}$.
Then, applying Lemmas \lemref{nzcomp} and \lemref{nz2},
we get the loop $Q_{r,s}$ defined by the table:

{ \newcommand\spx{\rule[-1pt]{0mm}{15pt}}  
\[
\begin{array}{c|c|c|c|c|}
 Q_{r,s}    &  z^j & a z^j  & b z^j  &  a b z^j \\
\hline
z^i  \spx  &  z^{i+j} & a z^{i+j} & bz^{i+j} &  ab z^{i+j} \\
\hline
a z^i \spx & a z^{i+j} &  z^{i+j+r} & abz^{i+j} &  b z^{i+j+r+2} \\
\hline
b z^i\spx  & b z^{i+j} & ab z^{i+j+1} &  z^{i+j+s} &  a z^{i+j+s+1} \\
\hline
a b z^i\spx  & ab z^{i+j} & b z^{i+j+r+1} & az^{i+j+s+2} &   z^{i+j+r+s +1} \\
\hline
\end{array}
\]
}

\noindent
Each $Q_{r,s}$ really defines a PACC-loop by Lemma \lemref{goodf2quot}.
$Q_{r,s}$ is never diassociative, since $aa\cdot b \ne a\cdot ab$.
There are $16$ possibilities for $r,s$, but up to isomorphism,
there are only two loops, $Q_{0,0}$  and $Q_{1,1}$:

If $r,s$ are both even, then $Q_{r,s} \cong Q_{0,0}$,
since if we let $2i = -r$ and $2j = -s$ and define 
$\hat a = a z^i$ and $\hat b = b z^j$,
then $(\hat a)^2 = (\hat b)^2 = 1$ and
$(\hat a\hat b) \ld (\hat b\hat a) = z$.
If we replace $a,b$ by $\hat a , \hat b$, we get the table for $Q_{0,0}$.

Also, if $r$ is even and $s$ is odd, then $Q_{r,s} \cong Q_{0,0}$.
To see this, let $2i = -r$ and $2j = -r-s-1$, and let
$\hat a = a z^i$ and $\hat b = ab z^j$.
Then $(\hat a)^2 = (\hat b)^2 = 1$ and
$(\hat a\hat b) \ld (\hat b\hat a) =
(a \cdot ab) \ld (ab \cdot a) = \hat z := z^{-1}$.
Replacing $a,b,z$ by $\hat a , \hat b, \hat z$, we again
get the table for $Q_{0,0}$.

Likewise, if $r$ is odd and $s$ is even, then $Q_{r,s} \cong Q_{0,0}$.

Finally, if $r,s$ are both odd, then $Q_{r,s} \cong Q_{1,1}$,
since we may let  $2i = -r + 1$ and $2j = -s + 1$ and define
$\hat a = a z^i$ and $\hat b = b z^j$; then
$(\hat a)^2 = (\hat b)^2 = z^1$, so that we get
the table for $Q_{1,1}$.

$Q_{0,0}$ and $Q_{1,1}$ are not isomorphic, since 
$\{x : x^2 = 1\}$ has size $2$ in $Q_{1,1}$ (namely, $\{1, z^2\}$),
and size $6$ in $Q_{0,0}$ (namely, $\{1, z^2, a, az^2, b, bz^2\}$).

The loops $Q_{r,s}$ are isomorphic to loop structures
defined on $\ZZ_4 \times \ZZ_2 \times \ZZ_2$ as follows:
\[
\begin{pmatrix}
x_1 \\ x_2 \\ x_3
\end{pmatrix}
\begin{pmatrix}
y_1 \\ y_2 \\ y_3
\end{pmatrix} :=
\begin{pmatrix}
x_1 + y_1 + r x_2 y_2 + s x_3 y_3 + 2 x_2 x_3 y_3 + 2 x_2 y_2 y_3 + x_3 y_2^2 \\
x_2 + y_2 \\
x_3 + y_3
\end{pmatrix} .
\]
The isomorphism is given on generators by
$z \leftrightarrow (1,0,0)^t$, $a \leftrightarrow (0,1,0)^t$, and
$b \leftrightarrow (0,0,1)^t$. The explicit formula was
found by making an \emph{ansatz} that the $\ZZ_4$-component
has the form $x_1 + y_1 + r x_2 y_2 + s x_3 y_3 + f(x_2,x_3,y_2,y_3)$,
and then computing the sixteen values of $f$ using the table.
Assuming further that $f$ is a homogeneous cubic polynomial
of the form
$f(x_2,x_3,y_2,y_3) = \sum_{0 \leq i\leq j\leq k\leq 1} \alpha_{ijk} x_i x_j y_k +
\sum_{0 \leq l\leq m\leq n\leq 1} \beta_{lmn} x_l y_m y_n$
leads to a system of nine linear equations with twelve unknowns in $\ZZ_2$,
namely the coefficients $\alpha_{ijk}$ and $\beta_{lmn}$.
The particular $f$ chosen here, namely
$f(x_1,x_2,y_1,y_2) = 2 x_2 x_3 y_3 + 2 x_2 y_2 y_3 + x_3 y_2^2$,
maximizes (though not uniquely) the number of zero coefficients.

Next, we consider Case (2).
Let $Z = \{1,c\}$ and $N = \{1,c,u,v\}$.
Let $Q = N\ \dot{\cup}\ Na\ \dot{\cup}\ Nb\ \dot{\cup}\ Nab$;
so $N \dot{\cup} Na$ is an 8-element group.  Since
$u,v \in N \backslash Z$, WLOG $N\ \dot{\cup}\ Na$ is nonabelian.
There are now three subcases:
\begin{itemize}
\item[2.1.]  $N \cong \ZZ_4$ and $N \dot \cup Na$ is the quaternion group.
\item[2.2.]  $N \cong \ZZ_4$ and $N \dot \cup Na$ is dihedral.
\item[2.3.]  $N$ is an elementary abelian $2$-group and $N \dot \cup Na$ is dihedral.
\end{itemize}
We shall see that Subcases (2.1) and (2.2) are impossible.
Note that in all three cases, $Z(N \dot \cup Na) = \{1,c\} = Z(Q)$,
and that all squares in $N \dot \cup Na$ lie in $\{1,c\}$.

By Lemma \lemref{index4}, $(a,a,b) = (a,b,b) = c$.
Using $\LCC$, we get
$a\cdot ba = [(ab) \rd a] \cdot a^2$, so
$(a\cdot ba) \cdot a = ab \cdot a^2$.
Let $ab = dba$, with $d \in N$.  Then, since $a^2, c$ are central,
\[
a^2 db = db a^2 = dba \cdot a \cdot c = ab \cdot a \cdot c =
a\cdot ba = a \cdot d\iv ab = a d\iv a\cdot b \cdot c \ \ ;
\]
the last `=' used the fact that $(a,d\iv b,b) =  (a,b,b)  = c$.
Thus, $a^2 d = a d \iv a c$, so that
$a d a\iv = d\iv c$, so $d$ is $u$ or $v$.

This now refutes Subcase (2.1), since in the quaternions (where
$c$ is now $-1$), the
conjugate of $i$ by any element other than $\pm 1, \pm i$
is $i\iv = -i$, not $- i\iv$.
It also refutes Subcase (2.2).
We shall view the dihedral group concretely as the symmetry
group of the square.
Here, $c$ is rotation by $180^\circ$,
$u$ and $v$ are rotations by $90^\circ$ and $270^\circ$,
and $a$ is a reflection, so the conjugate of $d \in \{u,v\}$ by $a$
is $d\iv$, not $d\iv c$.

\begin{table}[htb]
{ 
\footnotesize
\arraycolsep=2pt
\newcommand\spx{\rule[-6pt]{0mm}{19.5pt}}  
\[
\begin{array}{c|c|c|c|c|c|c|c|c|c|c|c|c|c|c|c|c|   }
 Q     & 1& c& u& v& a&ca&ua&va& b&cb&ub&vb&ab&cab&uab&vab \\
\hline
 1     & 1& c& u& v& a&ca&ua&va& b&cb&ub&vb&ab&cab&uab&vab \spx\\
\hline
 c     & c& 1& v& u&ca& a&va&ua&cb& b&vb&ub&cab& ab&vab&uab \spx\\
\hline
 u     & u& v& 1& c&ua&va& a&ca&ub&vb& b&cb&uab&vab& ab&cab \spx\\
\hline
 v     & v& u& c& 1&va&ua&ca& a&vb&ub&cb& b&vab&uab&cab& ab \spx\\
\hline
  a    & a&ca&va&ua& 1& c& v& u& ab&cab&vab&uab&cb& b& ub&vb  \spx\\
\hline
 ca    &ca& a&ua&va& c& 1& u& v&cab& ab&uab&vab& b&cb& vb&ub  \spx\\
\hline
 ua    &ua&va&ca& a& u& v& c& 1&uab&vab&cab& ab&vb&ub &  b&cb \spx\\
\hline
 va    &va&ua& a&ca& v& u& 1& c&vab&uab& ab&cab&ub&vb & cb& b \spx\\
\hline
  b    & b&cb&vb&ub&uab&vab&cab& ab& 1& c& v& u&ua&va&ca& a  \spx\\
\hline
 cb    &cb& b&ub&vb&vab&uab& ab&cab& c& 1& u& v&va&ua& a&ca  \spx\\
\hline
 ub    &ub&vb&cb& b& ab&cab&vab&uab& u& v& c& 1& a&ca&va&ua  \spx\\
\hline
 vb    &vb&ub& b&cb&cab& ab&uab&vab& v& u& 1& c&ca& a&ua&va  \spx\\
\hline
  ab   & ab&cab&uab&vab&vb&ub&cb& b&ca& a&va&ua& v& u& c& 1  \spx\\
\hline
 cab   &cab& ab&vab&uab&ub&vb& b&cb& a&ca&ua&va& u& v& 1& c  \spx\\
\hline
 uab   &uab&vab& ab&cab&cb& b&vb&ub&va&ua&ca& a& c& 1& v& u  \spx\\
\hline
 vab   &vab&uab&cab& ab& b&cb&ub&vb&ua&va& a&ca& 1& c& u& v  \spx\\
\hline
\end{array}
\]
}  
\caption{Subcase (2.3)}
\tablelabel{sixteen}
\end{table}

We are now left with Subcase (2.3).  Again,
$c$ is rotation by $180^\circ$, and WLOG,
$u$ is reflection on the $y$ axis,
$v$ is reflection on the $x$ axis,
and $a$ is reflection on the line $x=y$.
Then $a^2 = 1$, and we can assume WLOG that $d = u$, so $ab = uba$.
Furthermore, if $b,u$ commute then $ab, u$ do not commute.
Replacing $b$ by $ab$, we may assume WLOG, that $b,u$ do not
commute so that  $N \dot \cup Nb$ is nonabelian, and is thus
dihedral, as we have just seen.  Since some element of $Nb$ has
order $2$, we may assume WLOG that $b^2 = 1$.
In the dihedral group $N \dot \cup Nb$, all commutators are $1$ or $c$,
so that $ub = cbu$; it follows that $u$ (and also $v$) commute with $ab$.
We now can compute a complete table of $Q$; see
Table \tableref{sixteen}.  This table can mostly be filled out
using the commutation and association relations already described.
To fill out the lower right $4 \times 4$, we need to know
that $ab\cdot ab = v$.  To see that, use $\LCC$ to get
$ab\cdot ab = [(ab\cdot a)/(ab)]\cdot (ab\cdot b)$.
But $ab \cdot a = c \cdot a \cdot ba = c a u \cdot ab$,
and $ab\cdot b = ca$, so that $ab\cdot ab = cau \cdot ca = aua = v$.

Thus, Case (2) of Lemma \lemref{16cases} yields just one non-extra
PACC-loop of order $16$. Of course, one must verify that the loop
described by Table \tableref{sixteen} really is PACC.  Unlike
Case (1), this does not follow by the results of Section \secref{ext};
but it can easily be verified by a short computer program.

The loop of Case (2) is isomorphic to a loop structure defined on
$\ZZ_2\times \ZZ_2\times \ZZ_2\times \ZZ_2$ as follows:
\[
\begin{pmatrix}
x_1 \\ x_2 \\ x_3 \\ x_4
\end{pmatrix}
\begin{pmatrix}
y_1 \\ y_2 \\ y_3 \\ x_4
\end{pmatrix} :=
\begin{pmatrix}
x_1 + y_1 + x_3 y_2 + x_3 y_3 y_4 + x_4 y_2^2 + x_3 x_4 (y_3 + y_4)  \\
x_2 + y_2 + x_4 y_3^2 \\
x_3 + y_3 \\
x_4 + y_4
\end{pmatrix} .
\]
The isomorphism is given on generators by
$c \leftrightarrow (1,0,0,0)^t$, $u \leftrightarrow (0,1,0,0)^t$, 
$a \leftrightarrow (0,0,1,0)^t$, $b \leftrightarrow (0,0,0,1)^t$.
The term $x_3 y_2$ in the first component is determined by assuming
it is quadratic in the variables, and then using the equation
$au = va$ in the dihedral group $N \dot{\cup} Na$,
the upper left $8\times 8$ corner of the table. Analogously to Case (1),
the remaining terms in the first and second components were found
by assuming that they are homogeneous cubic polynomials in
$x_2, x_3, x_4, y_2, y_3, y_4$, where each term contains at
least one $x_4$ or $y_4$. Since $u\in N$, it is clear from the
table that the values of the polynomials are independent of $x_2$.
Using the table to compute values determines some coefficients;
the choice above maximizes (though not uniquely) the number of
zero coefficients.

\section{Order 27}
\seclabel{27}

\begin{table}[htb]
{ 
\tiny
\newcommand\mb[1]{\makebox[33pt]{$#1$}}  
\arraycolsep=1pt
\newcommand\spx{\rule[-14pt]{0mm}{33pt}}  
\newcommand\p[2]{\parbox{33pt}{\begin{center} $#1$\\$#2$ \end{center}}}  
\newcommand\z{\mathord{+}}  

\[
\begin{array}{c|c|c|c|c|c|c|c|c|c||c|   }
 Q   &  \mb{1}& \mb{a}& \mb{a^2}
    & \mb{b}& \mb{ab}& \mb{a^2b}
     &\mb{b^2}&\mb{ab^2}& \mb{a^2b^2 }&\  \mathrm{cubes}\  \\   
\hline
 1     &  0& 0 & 0 & 0 & 0 & 0 &0 &0 & 0 & 0 \spx\\
\hline
 a     &  0& 0 &\alpha& 0 & \gamma  & \p{\alpha}{+2\gamma}
    & 0&  2\gamma   & \p{\alpha}{+\gamma}  & \alpha \spx\\
\hline
 a^2   & 0&\alpha&\alpha &0&\p{\alpha}{+2\gamma} &  \p{\alpha}{+\gamma}
   & 0& \p{\alpha}{+\gamma} &\p{\alpha}{+2\gamma} & 2\alpha \spx\\
\hline
 b   & 0& \theta  &  \p{2\theta }{+\gamma}& 0&  \theta  & \p{2\theta }{+\gamma}
     & \beta &\theta \z \beta & \p{2\theta\z\beta}{+\gamma}& \beta  \spx\\
\hline
 ab  &0& \theta &\p{2\theta \z \alpha}{+ \gamma}&2\delta
       &\p{\theta }{+\gamma \z 2\delta }
 &\p{2\theta\z\alpha}{+ 2\delta}&\p{\beta}{+\delta}
        &\p{\theta\z\beta}{+2\gamma\z\delta}
   & \p{2\theta\z\alpha\z\beta}{+2\gamma\z\delta} & \alpha\z\beta \spx\\
\hline
   a^2b  &0& \theta\z\alpha & \p{2\theta\z\alpha}{+\gamma }
   & \delta &\p{\theta\z\alpha}{+2\gamma\z\delta}
       & \p{2\theta\z\alpha}{+2\gamma \z\delta}
    &\p{ \beta}{ + 2\delta} &\p{\theta\z\alpha\z\beta}{+\gamma\z2\delta }
    &\p{2\theta\z\alpha\z\beta}{+2\delta}  & 2\alpha\z\beta  \spx\\
\hline
 b^2   & 0& \p{2\theta }{+ 2\delta }&\p{\theta}{+ 2\gamma \z \delta }
    &\beta & \p{2\theta\z\beta}{+2\delta}  &\p{\theta\z\beta}{+2\gamma\z\delta}
    &\beta &\p{2\theta\z\beta}{+2\delta}
     &\p{\theta\z\beta}{+2\gamma\z\delta}&2\beta \spx\\
\hline
 ab^2  & 0& \p{2\theta}{+2\delta} &\p{\theta\z\alpha}{+2\gamma\z\delta}
    &\p{\beta}{+\delta} & \p{2\theta\z\beta}{+\gamma}
    & \p{\theta\z\alpha\z\beta}{+\gamma \z 2\delta}
   & \p{\beta}{+2\delta} &\p{2\theta\z\beta}{+2\gamma\z\delta}
     & \theta\z\alpha\z\beta 
   &  \alpha \z 2\beta \spx\\
\hline
 a^2b^2&0& \p{2\theta\z\alpha}{+2\delta}&\p{ \theta\z\alpha}{+2\gamma\z\delta }
  &\p{\beta}{+2\delta}&\p{2\theta\z\alpha\z\beta}{+2\gamma\z\delta}
    &\theta\z\alpha\z\beta 
      &\p{\beta}{+\delta}&\p{2\theta\z\alpha\z\beta}{+\gamma}
      &\p{\theta\z\alpha\z\beta}{+\gamma\z2\delta}  & 2\alpha\z2\beta  \spx\\
\hline
\end{array}
\]
}  
\caption{Order 27}
\tablelabel{27}
\end{table}

If $Q$ is a nonassociative PACC-loop of order 27,
then $|N(Q)| = |Z(Q)| = 3$, and $Q/N$ is an abelian group of exponent 3.
Thus $Q$ is an \emph{extraspecial} CC-loop, and one may then
in principle use Dr\'{a}pal's description of all extraspecial CC-loops
to classify the PACC-loops of order 27; see \S7 of \cite{DRAEX}.
Here we adopt a direct approach based on Lemmas \lemref{nzcomp2}
and \lemref{goodfquot}.

Say $N = \{1,n,n^2\}$ and $Q/N$ has generators $Na, Nb$.
Then $Q$ is determined by five parameters,
$\theta,\alpha,\beta,\gamma,\delta \in \ZZ_3 = \{0,1,2\}$,
where $ba = ab n^\theta$, $a^3 = n^\alpha$, $b^3 = n^\beta$,
$(b)E_a = b n^\gamma$, and $(a)E_b = a n^\delta$.
Then $a^i b^j \cdot a^k b^\ell = a^{i+k} b^{j + \ell} n^{f(i,j,k,\ell)}$,
and Table \tableref{27} displays $f(i,j,k,\ell)$.
Likewise, $(a^i b^j)^3 = n^{g(i,j)}$, and the table displays $g(i,j)$.
The table is computed using Lemma \lemref{nzcomp2},
and the loop is PACC by Lemma \lemref{goodfquot}.

We now count the number of distinct PACC-loops we have.

Let $T = \{x \in Q: x^3 = 1\}$.
By inspection of the table, we see that $T$ is always a subloop
of order either $9$ or $27$.

Next, let $M = M(Q) = \{x \in Q: E_x = I\}$, the set of Moufang elements of $Q$.  Clearly $N \subseteq M \subseteq Q$.
In fact, $N \subsetneqq M \subsetneqq Q$.
First, $M \ne Q$: $Q$ cannot be Moufang because it is not even WIP by Theorem \thmref{16or27} (or, by Chein \cite{CH}, a Moufang loop of order $27$ is a group).
Now, suppose that $M = N$.
Then $\gamma,\delta\in\{1,2\}$, so that $\gamma = \pm \delta$.
Referring to the table, we see:
\begin{align*}
& (a \cdot ab) \cdot ab  =
a^2 b \cdot ab \cdot n^\gamma  =
a^3 b^2 \cdot n^{\theta + \alpha + \delta} =
b^2 \cdot n^{\theta + 2\alpha + \delta}\ \  . \\ 
& a \cdot (ab)^2  =
a \cdot a^2 b^2 \cdot n^{\theta + \gamma + 2\delta} =
a^3 b^2 \cdot n^{\theta + \alpha + 2\gamma + 2\delta} =
b^2 \cdot n^{\theta + 2 \alpha + 2\gamma + 2\delta} \ \  .\\ 
& (a \cdot a^2b) \cdot a^2b  =
a^3 b \cdot a^2b \cdot n^{\alpha  + 2\gamma} =
b \cdot a^2b \cdot n^{2 \alpha  + 2\gamma} =
a^2b^2 \cdot n^{2 \theta + 2 \alpha} \ \  .\\
& a \cdot (a^2b)^2 =
a \cdot a^4b^2 \cdot n^{2\theta + \alpha + 2\gamma + \delta} =
a \cdot ab^2 \cdot n^{2\theta + 2 \alpha + 2\gamma + \delta} =
a^2 b^2 \cdot n^{2\theta + 2 \alpha + \gamma + \delta}\ \  .
\end{align*} 
If $\gamma = \delta$, then $(a \cdot ab) \cdot ab = a \cdot (ab)^2$, so $(a)E_{ab} = a$ (using $tx\cdot x = t x^2 \to t E_x = t$ by \RAlt = \RIP\ from Lemma \lemref{eqns}). 
Since also $(ab)E_{ab} = ab$, we have $E_{ab} = I$, contradicting $M = N$.  Likewise, if $\gamma = -\delta$, then $(a \cdot a^2b) \cdot a^2b = a \cdot (a^2b)^2$, which implies $E_{a^2b} = I$, again a contradiction.

As we noted in Conjecture \conjref{moufang}, we do not think that for arbitrary PACC-loops, the set of Moufang elements is necessarily a subloop. However, for this particular loop, $M$ is indeed a subloop of order $9$.
To see this, use $M \ne N$, and choose
a generator $a$ such that $E_a = I$ (that is, $\gamma = 0$).
Then $E_{a^2} = E_a^4 = I$, so $M \supseteq N \cup Na \cup Na^2$.
If $|M| > 9$, then we could choose the other generator $b$ so that $E_b = I$; but then $\delta = 0$, and then, as above, all $E_x = I$, so that $M = Q$, which is false.

Thus, $|M| = 9$, and as noted, we shall always choose our generators
so that $M = N \cup Na \cup Na^2$, and so $\gamma = 0$ and 
$\delta \in \{1,2\}$. Since $N = Z$, this shows that $M$ is a subloop.
We now have three cases:

\begin{itemize}
\item[I.]  $T = Q$.
\item[II.]  $T = M$.
\item[III.] $|T| = 9$ and  $T \ne M$.
\end{itemize}
Furthermore, each case will split into two subcases:
\begin{itemize}
\item[A.] $\exists x \in Q \backslash M \;  \exists y \in M \backslash N
\; [xy = yx]$.
\item[B.] $\forall x \in Q \backslash M \;  \forall y \in M \backslash N
\; [xy \ne yx]$.
\end{itemize}

In Case I, we have $\alpha = \beta = \gamma = 0$, and the generator
$b$ can be any element outside of $M$.
In Case IA, we can
choose $a$ and $b$ so
that $a,b$ commute, so that $\theta = 0$.  Then, WLOG $\delta = 1$,
since we can always replace $n$ by $n^2$.
In Case IB, we see from Table \tableref{27} that
$\theta \ne 0$ and
$\theta  + \delta \ne 0$  (in $\ZZ_3$), so that
$\theta = \delta \in \{1,2\}$. Again,
replacing $n$ by $n^2$, we may assume WLOG that $\theta = \delta = 1$.
Thus, Case I yields two loops.

In Case II, we have $T = M =  N \cup Na \cup Na^2$, so that
$\alpha = \gamma = 0$ and $\beta,\delta \in \{1,2\}$. 

In Case IIA, we can choose $a \in M \backslash N $ and
$b \in Q \backslash M$, so that $ab = ba$; so $\theta = 0$.
Now $b^3 = n^\beta$, but replacing $n$ by $n^2$ if necessary, WLOG $\beta = 1$.
Also, $(a)E_b = a n^\delta$, so 
$(a^2)E_b = a^2 n^{2\delta}$, so 
replacing $a$ by $a^2$ if necessary, WLOG $\delta = 1$;
note that by $\theta = \gamma = 0$,
$b a^2 = a^2 b$, so replacing $a$ by $a^2$ will not change our $\theta = 0$.

In Case IIB, note from the table that
$\theta \ne 0$ and $\theta + \delta \ne 0$ (in $\ZZ_3$).
WLOG $\theta = 1$, so that also $\delta = 1$;
that is, we choose $b,n$ so that $ba = abn$.
From the table, $b^2 a = a b^2 n$.  Also $(b^2)^3 = n^{2\beta}$.
Thus, replacing $b$ by $b^2$ if necessary, WLOG $\beta = 1$.

Thus, Case II yields two loops.

In Case III, we may choose generators so that
$M = N \cup Na \cup Na^2$ and $T = N \cup Nb \cup Nb^2$.
Then $\beta = \gamma = 0$ and $\alpha,\delta \in \{1,2\}$.

Now Case IIIA splits into two subcases:
\begin{itemize}
\item[1.] $\exists x \in T \backslash N \;  \exists y \in M \backslash N
\; [xy = yx]$.
\item[2.] $\forall x \in T \backslash N \;  \forall y \in M \backslash N
\; [xy \ne yx]$.
\end{itemize}
In case IIIA1, WLOG $\theta = 0$.
Then $b$ commutes with both of $a,a^2$, whereas $b^2$ commutes
with neither of $a,a^2$, so $b$ is fixed.
WLOG $\alpha = 1$; that is, we choose $a,n$ so that $a^3 = n$;
then also $a^2, n^2$ satisfy $(a^2)^3 = n^2$.
Now $\delta \in \{1,2\}$ is determined by $(a) E_b = (a)n^\delta$;
then $(a^2) E_b = (a^2)(n^2)^\delta$.  Thus, $\delta$ is determined
from the loop structure, so that we have two distinct loops.
In Case IIIA2, WLOG $\theta = 1$.
Then also $\delta = 1$, since otherwise $b^2$ would commute
with $a$.  But then regardless of the choice of $\alpha$,
neither of $a,a^2$ commute with any element of
$\{b, ab, a^2 b, b^2, ab^2, a^2 b^2\}$, contradicting the
assumption in IIIA.

In Case IIIB, WLOG $\theta = 1$; so $[b,a]$ is $n$, not $n^2$.
Then, as in Case IIIA2, we need $\delta = 1$, and this guarantees
that we are in IIIB for either choice of 
$\alpha \in \{1,2\}$.  
For each choice, we have
$\forall x \in M \backslash N\, \forall y \in T \backslash N\,
[\ [y,x] = x^{3\alpha}\ ]$,
so there are two distinct loops here.

Thus, Case III yields four loops, and there are exactly eight
nonassociative PACC-loops of order 27.
By Lemma \lemref{aip27}, the AIP holds in one if these loops iff
$\delta = \gamma + \theta$.
Since our loops all had $\gamma = 0$, 
we need $\delta = \theta$, which is true of the four ``B''
loops and false of the four ``A'' loops; so there are exactly
four nonassociative AIP PACC-loops of order 27.


\def\cyr{\fam\cyrfam\tencyr\cyracc}
\def\itcyr{\fam\cyrfam\tenitcyr\cyracc}


\begin{thebibliography}{99}

\bibitem{BAS2} {\cyr A. S. Basarab,
Ob odnom klasse} $G$-{\cyr lup},
{\itcyr Matematicheskie Issledovaniya}
{\cyr Tom} 3, {\cyr Vyp.} 2 (8)  (1968) 72--77.

\bibitem{BAS} {\cyr A. S. Basarab, Klass} LK-{\cyr lup},
{\itcyr Matematicheskie Issledovaniya} {\cyr Vyp.} 120 (1991) 3--7.
MR 92g:20115.

\bibitem{Bel} {\cyr V. D. Belousov},
{\itcyr Osnovy Teorii Kvazigrupp i Lup},
{\cyr Izdatel\cprime stvo <Nauka>, Moskva}, 1967.
MR 36{\#}1569, Zbl 163.01801.

\bibitem{Br} R. H. Bruck,
\textit{A Survey of Binary Systems},
Springer-Verlag, 1971.
MR 20{\#}76, Zbl 206.30301.

\bibitem{CH} O. Chein,
Moufang loops of small order, I,
\textit{Trans. Amer. Math. Soc.} \textbf{188} (1974) 31--51.
MR 48{\#}8673,  Zbl 286.20088.

\bibitem{CPS} 
O. Chein, H. O. Pflugfelder, and J. D. H. Smith, eds.,
\textit{Quasigroups and Loops: Theory and Applications}
Heldermann Verlag, 1990. 
MR 93g:20133,  Zbl 719.20036.

\bibitem{CD}
P. Cs\"{o}rg\H{o} and A. Dr\'{a}pal, 
Left conjugacy closed loops of nilpotency class $2$,
preprint.

\bibitem{DRA1} A. Dr\'{a}pal,
Conjugacy closed loops and their multiplication groups,
\textit{J. Algebra} \textbf{272} (2004),  838--850.

\bibitem{DRA2} A. Dr\'{a}pal,
Structural interactions of conjugacy closed loops,
preprint.

\bibitem{DRALCC} A. Dr\'{a}pal,
On multiplication groups of left conjugacy closed loops,
\textit{Comment. Math. Univ. Carolinae} \textbf{45} (2004) 191--212.

\bibitem{DRAEX} A. Dr\'{a}pal,
On extraspecial left conjugacy closed loops,
preprint.

\bibitem{FENA}  F. Fenyves, Extra loops I,
\textit{Publ. Math. Debrecen} \textbf{15} (1968) 235--238.
MR 38{\#}5976, Zbl 172.02401.

\bibitem{FENB}  F. Fenyves, Extra loops II,
\textit{Publ. Math. Debrecen} \textbf{16} (1969) 187--192.
MR 41{\#}7017, Zbl 221.20097.

\bibitem{GRA}  E. G. Goodaire and D. A. Robinson,
A class of loops which are isomorphic to all loop isotopes,
\textit{Canadian J. Math.} \textbf{34} (1982) 662--672.
MR 83k:20079, Zbl 467.20052.

\bibitem{GRB}  E. G. Goodaire and D. A. Robinson,
Some special conjugacy closed loops,
\textit{Canadian Math. Bull.} \textbf{33} (1990) 73--78.
MR 91a:20077, Zbl 661.20046.

\bibitem{KK} M. K. Kinyon and K. Kunen,
The structure of extra loops,
\textit{Quasigroups and Related Systems} \textbf{12} (2004) 39--60.

\bibitem{KKP} M. K. Kinyon, K. Kunen, and J. D. Phillips,
Diassociativity in conjugacy closed loops,
\textit{Comm. Algebra} \textbf{32} (2004) 767--786.

\bibitem{KUNC} K. Kunen,
The structure of conjugacy closed loops,
\textit{Trans. Amer. Math. Soc.} \textbf{352} (2000) 2889--2911.
MR 2000j:20132, Zbl 962.20048.

\bibitem{McM} W. W. McCune,
\textit{Mace\ 4.0 Reference Manual and Guide},
Argonne National Laboratory
Technical Memorandum ANL/MCS-TM-264, 2003;
\texttt{http://www.mcs.anl.gov/AR/mace4/}

\bibitem{Otter} W. W. McCune,
\textit{OTTER\ 3.3 Reference Manual and Guide},
Argonne National Laboratory Technical Memorandum ANL/MCS-TM-263, 2003;
{\tt{http://www.mcs.anl.gov/AR/otter/}}

\bibitem{NS} P. T. Nagy and K. Strambach,
Loops as invariant sections in groups, and their geometry,
\textit{Canad. J. Math.} \textbf{46} (1994)  1027--1056. 
MR 95h:20088, Zbl 814.20055.

\bibitem{Pf} H. O. Pflugfelder,
\textit{Quasigroups and Loops: Introduction},
Sigma Series in Pure Math. \textbf{8}, Heldermann Verlag, Berlin, 1990.
MR 93g:20132, Zbl 719.20036.

\bibitem{SO} {\cyr L. R. So\u\i kis, O Spetsial\cprime nykh lupakh},
in {\itcyr Voprosy Teorii Kvazigrupp i Lup}
({\cyr V. D. Belousov}, ed.),
{\cyr Redakc.-Izdat. Otdel Akad. Nauk Moldav. SSR, Kishinev}, 1970,
pp.~122--131.

\end{thebibliography}
\end{document}